\RequirePackage{fix-cm}

\documentclass[smallextended]{svjour3}       

\smartqed  

\usepackage{epsfig}
\usepackage{fancybox}
\usepackage{multicol}
\usepackage{color}
\usepackage{graphicx}
\usepackage{url}
\usepackage{floatrow}
\usepackage{subfig}
\usepackage{tikz}

\usepackage{booktabs}

\usepackage{amsmath}
\usepackage{amssymb}

\usepackage{amsthm}

\usepackage{lineno}

\makeatletter
\newcommand{\Rmnum}[1]{\expandafter\@slowromancap\romannumeral #1@}
\makeatother

\usepackage{multirow}

\usepackage{comment}

\begin{document}

\title{Numerical study on the effect of geometric approximation error in the numerical solution of PDEs using a high-order curvilinear mesh}

\titlerunning{Geometric error in numerical PDEs}        

\author{ Sehun Chun \and Julian Marcon \and Joaquim Peir{\'{o}} \and Spencer J. Sherwin}

\institute{ S. Chun \at Underwood International College, Yonsei University, South Korea, \\
               Tel.: +82-32-749-3607 \\
               \email{sehun.chun@yonsei.ac.kr}
               \and 
               J. Marcon \at 
               Aeronautics, Imperial College London, United Kingdom, \\
               \email{julian.marcon14@imperial.ac.uk} 
               \and
               J. Peir{\'{o}} \at 
               Aeronautics, Imperial College London, United Kingdom, \\
               \email{j.peiro@imperial.ac.uk} 
               \and 
               S. J. Sherwin \at
               Aeronautics, Imperial College London, United Kingdom, \\
               \email{s.sherwin@imperial.ac.uk} 
               }
               
\date{Received: date / Accepted: date}

\maketitle

\begin{abstract}
When time-dependent partial differential equations (PDEs) are solved numerically in a domain with curved boundary or on a curved surface, \textit{mesh error} and \textit{geometric approximation error} caused by the inaccurate location of vertices and other interior grid points, respectively, could be the main source of the inaccuracy and instability of the numerical solutions of PDEs. The role of these geometric errors in deteriorating the stability and particularly the conservation properties are largely unknown, which seems to necessitate very fine meshes especially to remove geometric approximation error. This paper aims to investigate the effect of geometric approximation error by using a high-order mesh with negligible geometric approximation error, even for high order polynomial of order $p$. To achieve this goal, the high-order mesh generator from CAD geometry called \textit{NekMesh} is adapted for surface mesh generation in comparison to traditional meshes with non-negligible geometric approximation error. Two types of numerical tests are considered. Firstly, the accuracy of differential operators is compared for various $p$ on a curved element of the sphere. Secondly, by applying the method of moving frames, four different time-dependent PDEs on the sphere are numerically solved to investigate the impact of geometric approximation error on the accuracy and conservation properties of high-order numerical schemes for PDEs on the sphere.
\end{abstract}

\section{Introduction}

The method of moving frames (MMF) was originally a continuous group theory developed by \'{E}lie Cartan to study the submanifolds of homogeneous spaces \cite{Cartan1,Cartan2,Cartan3}, and has been more developed and expanded in modern days for more practical and computational purposes, mostly in computer vision \cite{Olver1998,Olver2001,Faugeras1994} and medical sciences \cite{Piuze2015}. One of several efforts to adapt this framework to numerically solve partial differential equations on curved surfaces such as conservation laws \cite{MMF1}, diffusion equations \cite{MMF2}, shallow water equations \cite{MMF3}, and Maxwell's equations \cite{MMF4} was made by Chun.

By constructing an orthonormal unit vector, called \textit{moving frames}, at each grid point, the vector or the gradient of the scalar is expanded in these frames, which consequently yields the condition that requires corresponding equations to lie on the curved domain. This simple adaptation of moving frames removes non-Euclidean computational and geometrical redundancies, such as metric tensor and geometric singularity, that often deteriorate the accuracy and robustness of numerical schemes. In the modern era of computations involving complex geometry, many numerical schemes have been developed to solve differential equations on curved surfaces; however, in terms of accuracy, energy conservation, or stability, only the method of moving frames provides a unified and competitive framework for solving any partial differential equations on any type of curved surfaces that may exhibit anisotropy or angular rotation \cite{MMF1,MMF2,MMF3,MMF4}.

\begin{figure}[h]
\centering{
\includegraphics[
width=10cm]{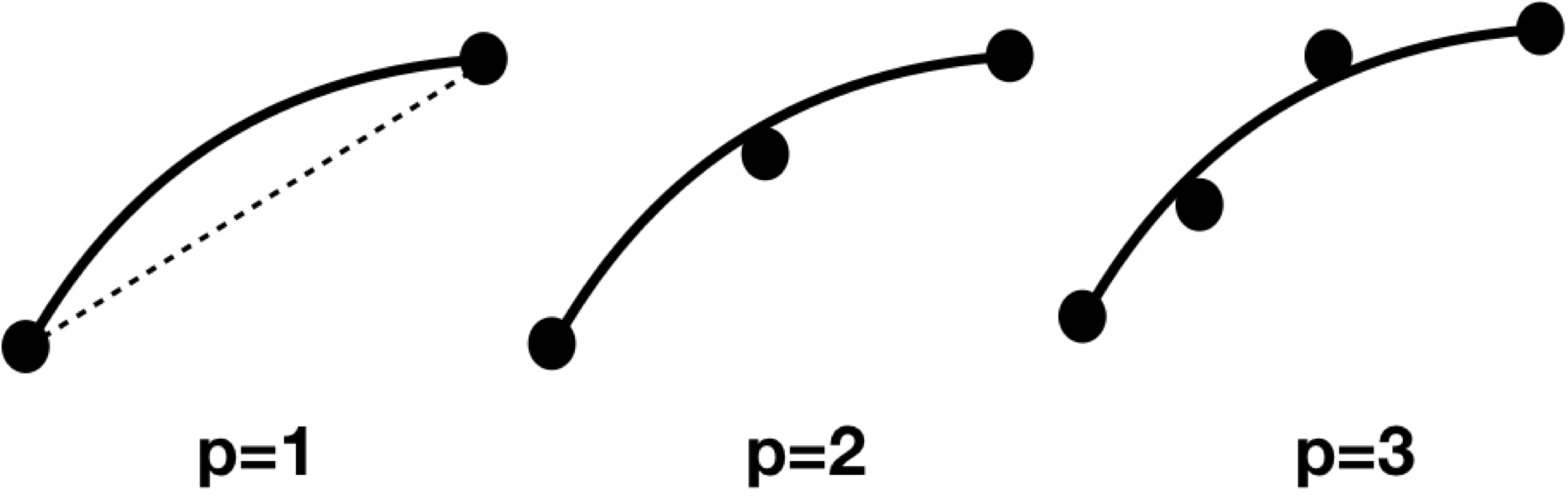}   }
\caption{Geometric approximation error: Interior grid point does not always lie on the curved domain, and therefore, the geometric approximation error does not decrease as $p$ increases.}
\label{MeshError}
\end{figure}

However, the success of moving frames in solving PDEs on curved surfaces was overshadowed by an \textit{inaccurate} curved mesh, which is referred as  a \textit{geometric approximation issue} for curvilinear mesh. In regular mesh of two-dimensional or three-dimensional Euclidean space, all the grid points lie in the Euclidean space irrespective of the polynomial order ($p$) used. On the other hand, some grid points in curvilinear meshes may not lie on the surface, particularly for those grid points other than vertices for $p \ge 3$. The location of vertices is fixed to remain on the surface through a mesh generator, of which inaccuracy corresponds to \textit{mesh error}. Additional grid points required for high-order polynomial approximation do not normally lie on the boundaries of the curvilinear element, of which inaccuracy corresponds to \textit{geometric approximation error}. Some meshing schemes do not fix the location of vertices, but relax it to minimize the maximum distance deviation of additional grid points from the surface. Nevertheless, in both cases, the geometric approximation issue is similar in the sense that the geometric approximation error does not decrease as the polynomial order increases. Figure \ref{MeshError} illustrates that the additional grid points for higher $p$, excluding the end points, are located outside the curved line as $p$ increases. Even when all the grid points lie on the surface, a question arises of whether the grid point distribution within the curved element is optimal for integration \cite{Hesthaven1}.

Numerical problems associated with geometric approximation error include the following: The first problem is that the geometric approximation error remains almost similar for all the $p$, and so the error in the numerical solution for high $p$ is mostly dominated by geometric approximation error and not by discretization errors. By considering that the method of moving frames displays the optimal accuracy of $p \ge 2$ \cite{MMF1}, non-converging geometric approximation errors with respect to $p$ frustrate the use of relatively high $p$ for finer resolution. The second problem involves the deterioration of numerical integration in a curved domain. Additional grid points lying inside or outside the surface are not identified as the best optimal nodal sets for the highest integration accuracy. This problem is closely related to the conservation properties of the corresponding numerical scheme, and therefore the numerical scheme on a curved mesh, even with a very small geometric approximation error, can result in excessive loss of conservative properties such as energy and mass after a long integration times, unless a very fine mesh is used.

The last problem is associated with the representation of anisotropy and angular rotation. The representation of anisotropy in the method of moving frames is denoted by the orientation and resizing of moving frames \cite{MMF2}, which should also be performed for surface rotation \cite{MMF3}. Moving frames are numerically constructed by differentiating the neighboring grid points. If the grid points do not lie on the surface, then the anisotropic tensor does not lie on the surface as a result. Consequently, the anisotropic direction and strength given at the location are inaccurately incorporated into the computational model because the anisotropy tensor should be projected on the surface along the direction of propagation.

To address these problems in this paper, the high-order surface mesh from CAD geometries, called \textit{Nekmesh} \cite{Peiro1,Peiro2,Turner}, is adapted to determine whether the aforementioned problems are solved on high-order mesh that  supposedly exhibits negligible geometric approximation error even for high $p$. Nekmesh is an open-source 3D high-order mesh generator within the \textit{Nektar++} framework \cite{Nektar++} that has been specifically utilised as a surface mesh generator in this study. This high-order mesh generator first deforms the CAD-based low-order linear mesh of the surface to conform to the geometry of the surface and reallocates the interior grid points through energy-based optimization for optimal nodal sets. Consequently, a curved mesh generated by Nekmesh satisfies the required conditions of high-order surface mesh as follows: For every $p$, (1) geometric approximation error decreases as $p$ increases, (2) grid point distribution is optimal for the integration accuracy.

In the context of the high-order mesh community, these tests can serve as an efficient \textit{quantification on the high-order mesh quality} because all the meshes of a sphere are supposed to be curvilinear meshes. All codes for the MMF schemes on Nekmesh are publicly available in Nektar++ version 4.5.0 or higher downloadable on GitHub\footnote{\texttt{https://gitlab.nektar.info} \cite{Nektar++}}. \texttt{Nekmesh} command in \texttt{build$/$utilities} folder of Nektar++ converts an mcf file into an xml file that can be used as a curvilinear \textit{mesh} of the sphere.

This paper is organized as follows. A brief summary of a high-order mesh generator, Nekmesh, is presented in Section 2. The high-order mesh properties of Nekmesh in comparison to other meshes are presented in Section 3. Section 4 explains the use of moving frames in the numerical solution of PDEs. In Section 5, the method of moving frames is applied to solve four differential operators, including divergence, gradient, and curl, on a curved element of a sphere to compare the accuracy of covariant differentiation on a high-order curved element of a sphere with that on traditionally available curved elements using projected mesh \cite{risser} and Gmsh \cite{Gmsh}. In Section 6, four partial differential equations such as conservation laws, diffusion equations, the shallow water equations, and Maxwell's equations, are solved on the high-order curved mesh of the sphere for analysis and comparison. Moreover, the effect of geometric approximation errors associated with conservation properties such as mass and energy especially are investigated for shallow water equations. The discussion is presented in Section 7.

\section{High-order mesh generator, \textit{NekMesh}}\label{sec:NekMesh}

\textit{NekMesh} is a set of tools to generate and manipulate high-order curvilinear meshes~\cite{Peiro2,Turner} as a part of the open-source \textit{Nektar++}~\cite{Nektar++} platform, a framework for CFD solvers based on the spectral/\textit{hp} element methods. \textit{NekMesh} was initially developed for converting mesh file formats but has evolved into an extensive set of tools for generating~\cite{Peiro2}, optimizing~\cite{Turner2017}, and adapting~\cite{Marcon2017} high-order meshes involving strict requirements for geometrical accuracy.

The generation of high-order meshes in NekMesh is based on an \textit{a posteriori} approach, which follows a bottom-up procedure proposed in~\cite{Peiro1}. A linear mesh is first obtained through traditional methods well-known by the low-order community. To achieve the desired high-order polynomial discretization of a high-order mesh, additional points are added to the linear mesh. These points are added along the edges and across faces, and then projected onto the boundaries to obtain a geometrically accurate boundary discretization in the CAD domain.

It is important that this discretized representation of the boundaries is accurate for machine precision because even small geometrical inaccuracy in the mesh would result in considerable loss of solver accuracy. These high-order meshes are therefore very sensitive to the underlying CAD boundary representation, which can lead to distorted boundary elements where the curvature of CAD boundaries is high. To solve this problem, {NekMesh} relies on an optimization procedure for the high-order nodes that lie on the boundaries to address this problem.

This procedure also follows a bottom-up approach where curve-bound edge interior nodes are first optimized, followed by surface-bound edge interior nodes, and finally face interior nodes. At each of these steps, the optimization of the location of high-order nodes relies on the minimization of the deformation energy of a virtual system of springs. In the first step, only edges lying on a CAD curve are processed and end nodes that correspond to the linear mesh nodes are fixed. The high-order interior nodes are then assembled in a spring system in parametric space and the deformation energy of the entire system is minimized.

The same procedure is applied in the second step where the spring system is built in the two-dimensional parametric space of the CAD surface on which the edge lies. End nodes are again fixed resulting in high-order edge interior nodes lying approximately on the geodesic between the two end points on the three-dimensional surface. Finally, face interior nodes are processed. We fix the nodes lying on the edges of the boundary elements and build a two-dimensional spring system in the parametric space of the surface. Nodes are connected to each of their neighbors; in the case of a triangle, each node is effectively connected to its six neighboring nodes arranged in a hexagon around it. The optimization of all high-order boundary nodes ensures that the surface mesh is smooth as long as the CAD representation is of good quality.

This procedure relies on a robust CAD system for querying the CAD model to obtain the location of projections onto surfaces or gradients required for the optimization procedure described above. In order to interface any CAD system, NekMesh was developed with a wrapper that allows developers to interact with CAD systems on a high level through a minimum number of functions calls. The default CAD system used by NekMesh is OpenCASCADE~\cite{OpenCascadeSAS2018}, but other systems have also been interfaced, including CADfix through its CFI interface for complex geometries~\cite{Turner2017b,Marcon2018}.

\section{Comparisons of mesh properties }

For the purpose of comparison, two different types of non high-order curved meshes are used: the first mesh is the projected mesh, called \textit{ProjMesh} \cite{risser}. The edges of a cube of unit length are equidistantly dissected according to a user-defined edge length, and the vertices of the derived edges are projected on a sphere. The second mesh is obtained by the built-in `sphere.geo' in \textit{Gmsh} \cite{Gmsh}, an open-source three-dimensional finite element mesh generator. For Gmsh, each edge is curved with a second order polynomial and the Frontal algorithm is used as the mesh algorithm.

We compare these three methods with respect to their features. First, ProjMesh projects the mesh of a cube onto an analytic sphere. This projection is exact but the relevance for industrial applications is limited, though it is one of the most popular approach for generating meshing on the sphere. Second, Gmsh permits the generation of meshes on arbitrary manifold by using a CAD system for the placement and the projection of nodes and vertices. Finally, NekMesh proposes an approach that is mindful of high-order information. NekMesh not only projects high-order nodes onto CAD surfaces, but also optimizes their location on geodesics, as explained in Sect.\ref{sec:NekMesh}. See Table \ref{table::CompMeshes} for quick comparison of those meshes. In the remaining section, the importance of this last step in NekMesh will be explained.

\begin{table}[ht]
\resizebox{\textwidth}{!}{
\begin{tabular}{c c c c }
 \hline\noalign{\smallskip} 
  & ProjMesh & Gmsh & NekMesh  \\ 
    \noalign{\smallskip}\hline\noalign{\smallskip}
  Type & projection & lower-order mesh & high-order mesh  \\
  Optimal grid distribution & No & Yes & Yes      \\
  Mesh error & Negligible & Negligible  & Negligible \\
 Geometric approximation  &  \multirow{ 2}{*}{Non-negligible} &  \multirow{ 2}{*}{Non-negligible} &  \multirow{ 2}{*}{Negligible}  \\
  error for $p \ge 2$ & & & \\ 
  \noalign{\smallskip}\hline\noalign{\smallskip}
\end{tabular}
\caption{ Comparison of curvilinear meshes used in the paper. \textit{Negligible} means that the corresponding quantity is relatively small in comparison to other discretization errors. }
\label{table::CompMeshes}
}
\end{table}

In all approaches, high-order nodes are added to linear edges and faces and then projected onto the manifold, either analytically or onto a CAD model. The distribution of these nodes, commonly uniform before projection, is not preserved upon projection. The projection effectively introduces an additional mapping between the linear and the curvilinear edges. The \textit{energy}, or also called the \textit{degree of deformation} of this mapping, reduces the convergence properties of the polynomial discretisation. The concept of energetic mappings is illustrated in Fig.~\ref{fig:mapping} where a high energy or highly deformed mapping (left) is compared to a low energy or low deformation mapping (right). The terminology \textit{energy} is introduced in analogy to Fourier transforms where a high energy mapping would be expected to have high frequency components in the mapping if the mapping is decomposed into hierarchical/Fourier expansions.

When the projection of points that are originally equispaced under the action of a mapping produce a distribution of points that remain equally spaced after projection then the energy content of the mapping is such that the energy of higher frequencies of the mapping is significantly lower than those of lower frequencies. By introducing a high energy mapping, this energy spectrum of the polynomial mapping is shifted to higher frequencies and the accuracy of the representation decreases.

In the example shown in Fig.~\ref{fig:mapping}, the points being equispaced before projection in the region ($0 \leq \xi \leq 1$) are projected to relatively uniformly spaced points by the low energy mapping $x^b(\xi)$ and a non-uniform distribution of points by the high energy mapping $x^a(\xi)$. If two new points are introduced as highlighted in red in Fig.~\ref{fig:mapping}, the low energy mapping $x^b(\xi)$ maintains the similar distribution when projected, as shown on the right plot where the high energy mapping $x^a(\xi)$ distributes these points in a more distorted manner after projection. A direct consequence of a high energy mapping is that the metric properties used in integrating and differentiating the spectral/hp element approximation need to be approximated with a high order polynomial expansion which are likely to lead to a degradation in the solution approximation when compared to a low energy mapping.

\begin{figure}
    \centering
    \begin{tikzpicture}
  \begin{scope}
    \draw (0,0) -- (4,0)  ;
    \fill[black] (0,0) circle (2pt) node[above] {0} ;
    \fill[black] (1,0) circle (2pt) ;
    \fill[black] (2,0) circle (2pt) ;
    \fill[black] (3,0) circle (2pt) ;
    \fill[black] (4,0) circle (2pt) node[above] {1} ;
    
    \fill[red] (1.2,0) circle (2pt) ;
    \fill[red] (3.2,0) circle (2pt) ;
    
    \draw[->] (0,-0.5) -- (1,-0.5) node[right] {$\xi$} ;
  \end{scope}
  
  \begin{scope}[xshift=-100,yshift=100]
    \draw (0,0) to [bend left] coordinate[pos=0.1] (A1) coordinate[pos=0.18] (B1) coordinate[pos=0.5] (A2) coordinate[pos=0.9] (A3) coordinate[pos=0.92] (B2) (4,2) ;
    \fill[black] (0,0) circle (2pt) node[left] {$x^a_0$} ;
    \fill[black] (A1) circle (2pt) ;
    \fill[black] (A2) circle (2pt) ;
    \fill[black] (A3) circle (2pt) ;
    \fill[black] (4,2) circle (2pt) node[above] {$x^a_1$} ;
    
    \fill[red] (B1) circle (2pt) ;
    \fill[red] (B2) circle (2pt) ;
  \end{scope}
  
  \begin{scope}[xshift=100,yshift=100]
    \draw (0,0) to [bend left] coordinate[pos=0.25] (A1) coordinate[pos=0.3] (B1) coordinate[pos=0.5] (A2) coordinate[pos=0.75] (A3) coordinate[pos=0.8] (B2) (4,2) ;
    \fill[black] (0,0) circle (2pt) node[left] {$x^b_0$} ;
    \fill[black] (A1) circle (2pt) ;
    \fill[black] (A2) circle (2pt) ;
    \fill[black] (A3) circle (2pt) ;
    \fill[black] (4,2) circle (2pt) node[above] {$x^b_1$} ;
    
    \fill[red] (B1) circle (2pt) ;
    \fill[red] (B2) circle (2pt) ;
  \end{scope}
  
  \draw[very thick,-latex] (-1,1) to[bend left] node[right] {$x^a \left( \xi \right)$} (-2,3) ;
  \draw[very thick,-latex] (5,1) to[bend right] node[left] {$x^b \left( \xi \right)$} (6,3) ;
  
\end{tikzpicture}
    \caption{Mapping from reference to physical space.}
    \label{fig:mapping}
\end{figure}

From geometrical considerations alone, it can be deduced that, the larger the local curvature of the element, the greater the energy of the mapping. The importance for NekMesh to optimise the distribution of high-order nodes after projection is now clearer. Through a spring system as shown in Sect.~\ref{sec:NekMesh}, NekMesh effectively redistributes nodes according to their original, and therefore optimal, distribution, which leads to removing the in-plane mapping error introduced by the projection onto the sphere. For example, consider $x^b \left( \xi \right)$ the lower energy result of the optimisation of mapping $x^a \left( \xi \right)$ in Fig.~\ref{fig:mapping}. The accuracy of the high-order mesh in NekMesh is further improved by choosing a Gauss-Lobatto-Legendre point distribution. This distribution has been shown to have a lower Lebesque constant, or $L_\infty$-norm approximation error, than an evenly-spaced node distribution, as used in ProjMesh before projection. One of the key features of NekMesh is best explained in terms of geometric approximation error that will be explained in the following subsection.

\begin{figure}[ht]
\centering
\subfloat[ProjMesh] {\label{Mesherr1} \includegraphics[
width=3.5cm]{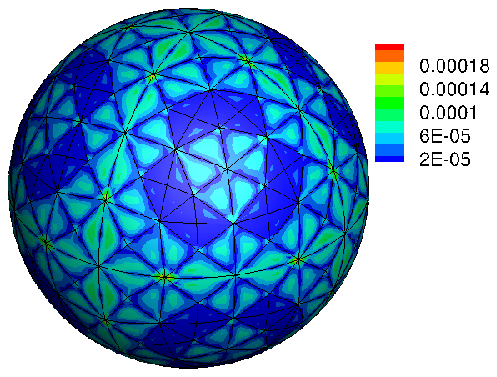} }
\subfloat[Gmsh] {\label{Mesherr2} \includegraphics[
width=3.5cm]{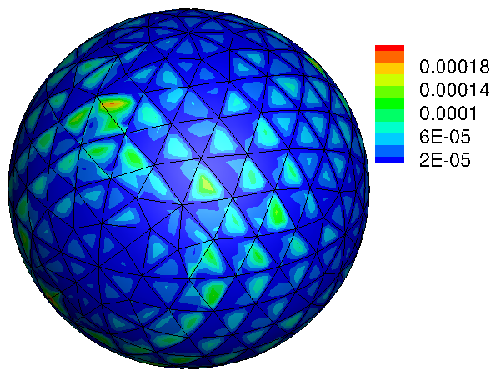} }
\subfloat[NekMesh] {\label{Mesherr3} \includegraphics[
width=3.5cm]{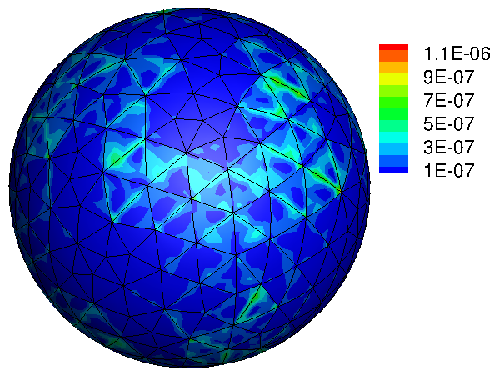} }
\caption{Distribution of geometric approximation error for a sphere. $h=0.4$, $p=4$.}
\label {MeshErrorSphere}
\end{figure}

\subsection{Geometric approximation error}

For the maximum edge length  $h$ for the sphere with radius of unit length, several curved meshes are obtained from separate algorithm. We define the mesh error in $L_2$ associated with the sphere as follows:
\begin{equation}
\mbox{Mesh error}  \equiv   \sqrt{ \sum_{v_i=1}^{N_v} \dfrac{( r_{\textrm{exact}} - {r}_{v_i} )^2}{N_v} }, \label{Mesherror}
\end{equation}
where ${r}_{v_i}$ is the radius of the $v_i$th \textit{vertex}, $r_{exact}$ is the exact radius of the sphere, and $N_v$ is the total number of vertices. On the other hand, geometric approximation error is defined similarly as follows. 
\begin{equation}
\mbox{Geometric approximation error}  \equiv    \sqrt{  \sum_{i=1}^N  \dfrac{( r_{\textrm{exact}} - r_i)^2}{N}   },~~~  \label{GAerror}
\end{equation}
where $r_i$ is the radius of the $i$th \textit{grid points} and $N$ is the total number of grid points. Note that mesh error is the same as geometric approximation error if all the grid points are vertices with no interior grid points. The distribution of geometric approximation error for the three meshes are displayed in Fig. \ref{MeshErrorSphere}.

\begin{figure}[ht]
\centering
\subfloat[versus $h$] {\label{gaeh} \includegraphics[
width=5cm]{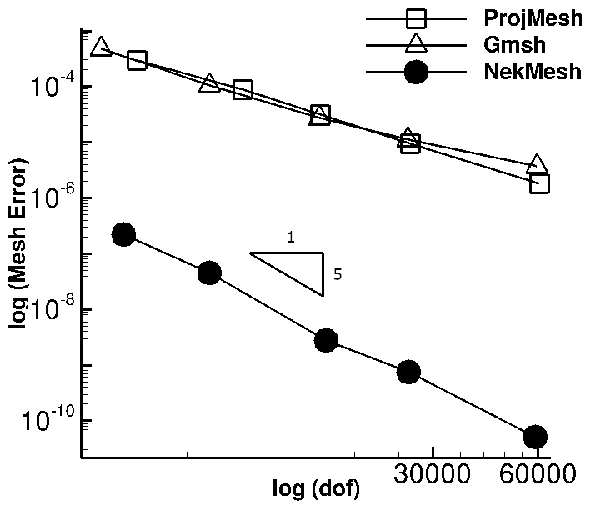} }
\subfloat[versus $p$] {\label{gaep} \includegraphics[
width=5cm]{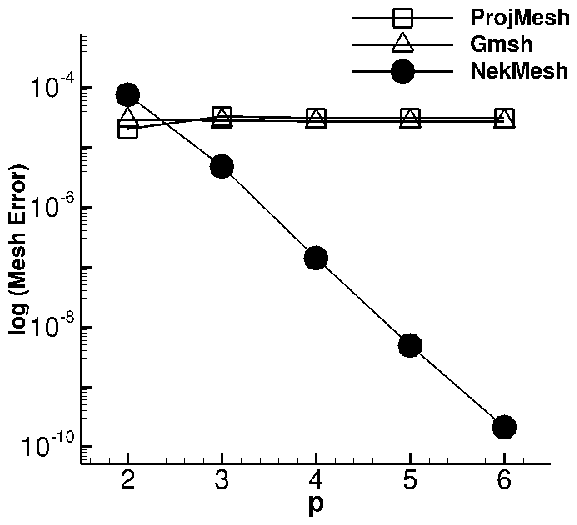} }
\caption{(a) Geometric approximation error versus edge length ($h$) given $p=4$ and (b) geometric approximation error versus $p$ given $h=0.4$.  dof = degree of freedom = number of total grid points. $4\times10^3$ corresponds to $h \approx 0.6$ and $6 \times 10^4$ corresponds to $h \approx 0.2$}
\label {MeshErrorph}
\end{figure}

Fig. \ref{MeshErrorph} and Table \ref{table::Mesherror} present the degree of freedom versus the geometric approximation error for three different meshes. Fig. \ref{gaeh} demonstrates that the geometric approximation error of ProjMesh and Gmsh is approximately the same, but the geometric approximation error of Nekmesh is significantly smaller than those of ProjMesh and Gmsh when $p=4$. This difference is slightly magnified when $h$ becomes smaller, i.e., the number of degrees of freedom (dof) is larger. When the dof is close to $4\times10^3$ ($h \approx 0.6$), the geometric approximation error in Nekmesh is approximately $7.5\times10^{-4}$ times smaller than the other two meshes' geometric approximation errors. When the dof is close to $6 \times 10^4$ ($h \approx 0.2$), the geometric approximation error in Nekmesh is approximately $2.8 \times 10^{-5}$ times smaller than the other two meshes's geometric approximation errors.

More perspectives can be obtained when the geometric approximation error is plotted against $p$. Fig. \ref{gaep} presents the convergence of geometric approximation error versus $p$ in $L_2$, which indicates that the vertices of the simplices all lie on the surface, but some of additional grid points added to the elements do not lie on the surface. Gmsh exhibits similar behavior, but the increased rate is much smaller. On the other hand, the geometric approximation error of NekMesh decreases as $p$ increases. Consequently, the $L_2$ error of ProjMesh and Gmsh remains approximately the same regardless of $p$, but the $L_2$ error of NekMesh decreases exponentially for $p \ge 2$. At $p=2$, Nekmesh exhibits the largest geometric approximation error, approximately $3.5$ times larger than that of ProjMesh. However, the geometric approximation error in Nekmesh decreases exponentially as $p$ increases. For example, at $p=6$, the geometric approximation error of NekMesh is $5.0 \times 10^{-5}$ smaller than those of the other two meshes.

\begin{table}[ht]
\begin{tabular}{c c c c c c c c}
 \hline\noalign{\smallskip} 
 $h$ &  \multicolumn{3}{c}{$p$-error} & $p$ & \multicolumn{3}{c}{$h$-error}   \\
 \cmidrule(lr){2-4} \cmidrule(lr){6-8}  
 & ProjMesh & Gmsh & NekMesh &  & ProjMesh & Gmsh & NekMesh \\
  \noalign{\smallskip}\hline\noalign{\smallskip}
 0.6 & 2.97E-04  & 4.81E-04 & 2.22E-07 & 2 & 2.05E-05 & 2.84E-05 & 7.48E-05 \\
 0.5 & 8.81E-05 & 1.05E-04 & 4.59E-08 & 3 & 3.32E-05 & 2.76E-05 & 4.77E-06  \\
 0.4 & 3.09E-05 & 2.74E-05 & 2.77E-09 & 4 & 3.09E-05 & 2.74E-05 & 1.39E-07  \\
 0.3 & 9.25E-06 & 1.11E-05 & 7.61E-10 & 5 & 3.06E-05 & 2.74E-05 & 4.94E-09  \\
 0.2 & 1.83E-06 & 3.71E-06 & 5.12E-11 & 6 & 3.04E-05 & 2.74E-05 & 2.10E-10 
\end{tabular}
\caption{Geometric approximation error in $L_2$ norm. $p$-error for $p$=$4$. $h$-error for $h$=$0.4$. $2.97E$-$04$ is shorthand for $2.97 \times 10^{-4}$.}
\label{table::Mesherror}
\end{table}

\section{Moving frames for PDEs on the sphere}

\begin{figure}[ht]
\centering
 \includegraphics[
width=4cm]{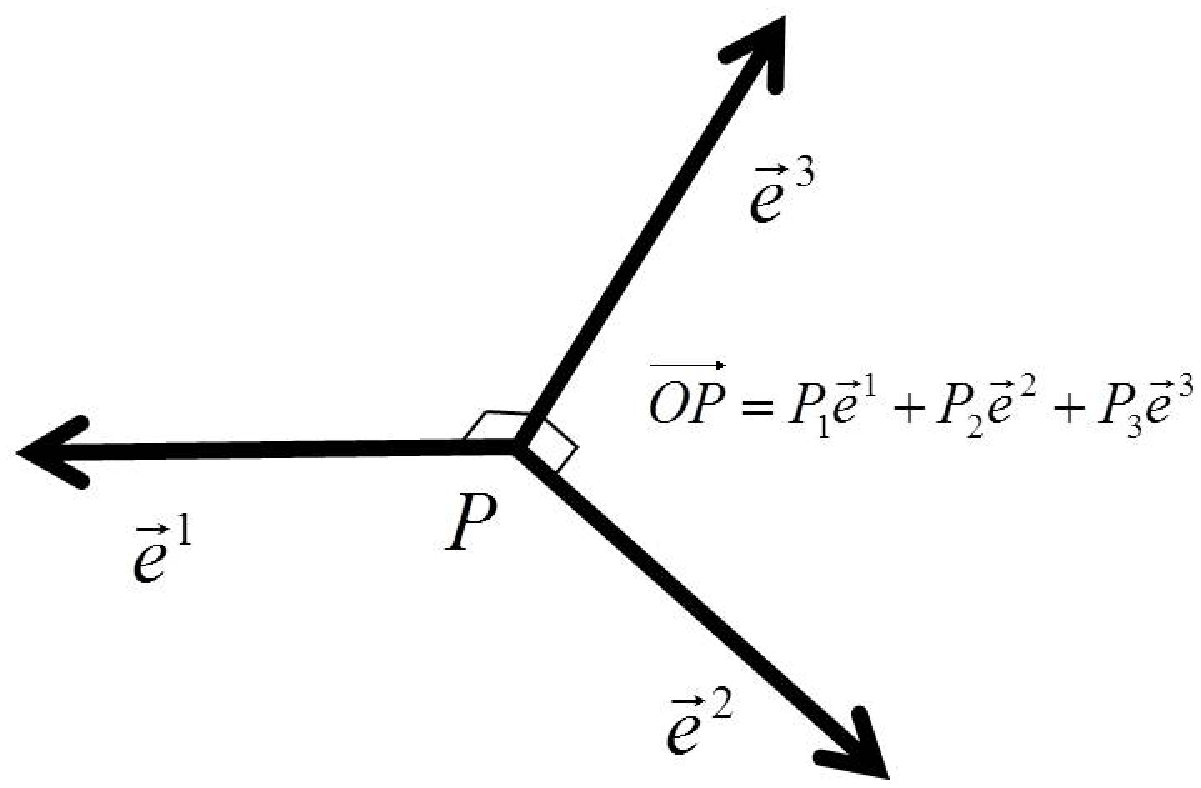}
\includegraphics[
width=3cm]{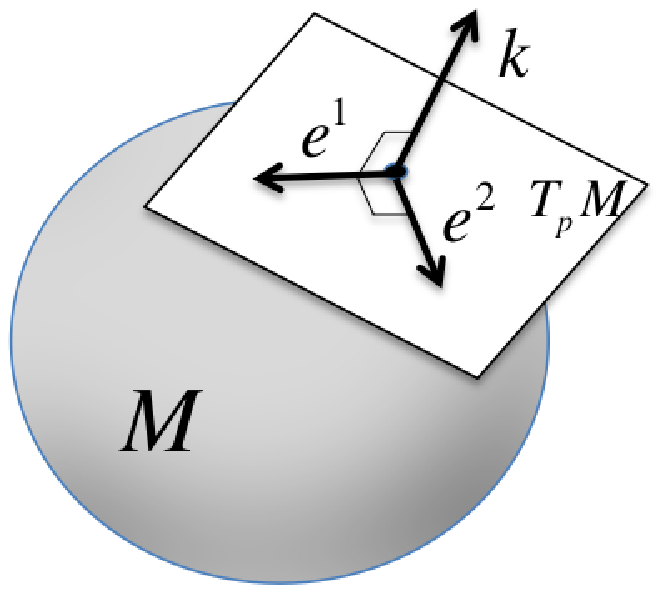}
\caption{Illustration of moving frames.}
\label {MFillust}
\end{figure}

Let $\mathbf{e}^i,~ 1 \le i \le 3$ be the moving frames that are constructed at each grid point $P$ as shown in Fig. \ref{MFillust} (the index that indicates the dependency of moving frames on grid points is omitted). Moving frames can be considered orthonormal bases, i.e., orthogonal to each other and of unit length such that
\begin{equation*}
\| \mathbf{e}^i \| = 1,~~~~~ \mathbf{e}^i \cdot \mathbf{e}^j = \delta^i_j, ~~~~1 \le i,j \le 3.
\end{equation*}
where $\delta^i_j$ is the Knocker delta.

An additional condition required for moving frames in the numerical solutions of PDEs includes the differentiability of $\mathbf{e}^i$ for each curved element. However, the differentiability of $\mathbf{e}^i$ for the entire domain is not required. The numerical construction of moving frames on curved surfaces are simple and described in details in prior studies \cite{MMF1,MMF2,MMF3}. For example, on a curved element, two moving frames $\mathbf{e}^1$ and $\mathbf{e}^2$ consequently lie on the tangent plane, whereas $\mathbf{e}^3$ is in the same direction of the surface normal vector $\mathbf{k}$, as shown in Fig. \ref{MFillust}. 

Adapting moving frames in the numerical solution of PDEs on a curved surface can be easily achieved by expanding the vectors or the gradient with $\mathbf{e}^i$. For example, let $\mathbf{v}$ be the velocity vector in the equations and be expanded by the moving frames at $P$ as follows
\begin{equation}
\mathbf{v} = v^1 \mathbf{e}^1 + v^2 \mathbf{e}^2,~~~~~~ v^i \in \mathbf{R}.   \label{vinMF}
\end{equation}

Let $\boldsymbol{\nu}^i$ be the tangent vector aligned along $\mathbf{e}^i$ such that $\boldsymbol{\nu}^i =  \beta \mathbf{e}^i$ where $\beta \neq 1$. The vector $\mathbf{v}$ can also be expanded in the tangent vectors $\boldsymbol{\nu}^i$ with a modified coefficient $v^i$, but maintaining the magnitude of the tangent vector in the basis vector often introduces the metric tensor in the PDEs. The direction of $\boldsymbol{\nu}^i$ cannot be chosen randomly because the length of $\boldsymbol{\nu}^i$, which can be easily obtained analytically from the known axis, should be computed with the metric tensor if the axis is not known, which is known to be computationally challenging. Using moving frames of unit length removes the metric tensor in numerical schemes, and therefore the direction can be chosen at random.

Representing vector or the gradient in moving frames does not produce additional numerical errors, but its derivative does. Let $\xi^i$ be the curved axis with the tangent vector $\mathbf{e}^i$. The divergence of $\mathbf{v}$ in Eq. \eqref{vinMF} is then denoted as
\begin{equation}
 \sum_{i=1}^2 \left [ \frac{\partial v^i}{\partial \xi^i} + \sum_{j=1}^2 \Gamma^i_{ij} v^j \right ] =  \nabla \cdot \mathbf{v}  = \sum_{i=1}^2 \left [ \nabla v^i \cdot \mathbf{e}^i + v^i \nabla \cdot \mathbf{e}^i \right ]  .  \label{divinMF}
\end{equation}

The left equality is the exact expression of the divergence in the curved domain and the right equality is the alternative derivation of divergence in moving frames. The first component in the right-hand side is the same because it is a metric-free component, but the second component is not the same because the divergence in the right-hand side is computed on the tangent plane as an Euclidean space.

In the numerical schemes with moving frames, $\nabla \cdot \mathbf{e}^i$ is computed on the tangent plane, not on curved surfaces, without using the metric tensor. Consequently, additional errors are generated in addition to the discretization error for computing divergence. This error is referred to as the \textit{MMF error}, which is associated with the error arising out of the curvature of the domain. However, it has been proved analytically and validated computationally that the MMF error is always negligible in comparison to the discretization error if all the curved elements of a curved domain have almost constant curvature \cite{MMF1,MMF2}. Because any modern mesh generator can produce a mesh that can adaptively tesselate any more highly curved region into smaller elements, the condition for the negligible MMF error can be easily achieved by a smaller $h$ for a region of higher curvature. In other words, if we increase $p$ or decrease $h$ to decrease the discretization error, the MMF error then decreases more rapidly and becomes negligible at a sufficient resolution. This feature has more advantages than the covariant formulation where the metric tensor is presented as an approximated coefficient, rather than the exact coefficient, that deteriorates the accuracy of the system.

\subsection{Numerical schemes}

For the test problems in the remainder of this paper, the weak formulation is adapted in the context of discontinuous Galerkin (DG) method. Moving frames is not restricted to the DG method, but the DG method is chosen according to the authors' preference. Any numerical scheme of high-order should yield the similar results. A brief description of the DG method is shown as follows: consider a conservation law for a variable $u$ and the velocity vector $\mathbf{v}$ as follows.
\begin{equation*}
\frac{\partial u}{\partial t} + \nabla \cdot \mathbf{v} = 0 .
\end{equation*}
The weak form of the above equation is obtained by multiplying a test function $\varphi$ and integrating over a sufficiently smooth domain $\Omega$ such as 
\begin{equation}
\int_{\Omega} \frac{\partial u}{\partial t} \varphi dx + \int_{\Omega} \nabla \cdot \mathbf{v} \varphi dx= 0 .  \label{DG1}
\end{equation}
Suppose $\Omega_h$ is a tessellation of the domain $\Omega$ into $N$ elements $\Omega^i_h; i=1, \ldots, N$, with characteristic edge length $h$ such as
\begin{equation*}    
\Omega _h  = \bigcup_i^N  \Omega^i_h,   \quad
\Omega^i_h \bigcap  \Omega^j_h = \emptyset~~\mbox{if}~ \; i \neq j
\end{equation*}
Moreover, consider the finite dimensional space consisting of discontinuous piecewise polynomial functions over $\Omega_h$ such that
\begin{equation*}
\mathcal{S}^p_h  = \{ u^h \in P^p (\Omega^i_h ), ~~~ \Omega^i_h \in \Omega_h \}
\end{equation*} 
Choose a basis function $\Phi_n \in \mathcal{S}^p_h$. Then, the solution and the test function are expressed as the linear combination of basis functions as follows.
\begin{equation}
u^h = \sum_{n} \hat{u}^h_n \Phi_n,~~~~ \varphi^h = \sum_{n} \hat{\varphi}^h_n \Phi_n .    \label{DG2}
\end{equation} 
By substituting Eqs. \eqref{DG2} into Eq. \eqref{DG1} and integration by parts, the conservation laws can be rewritten as
\begin{equation*}
\int_{\Omega} \frac{\partial \hat{u}_h}{\partial t} \hat{\varphi}_h dx - \int_{\Omega} \nabla \hat{\varphi}_h  \cdot \mathbf{v}_h dx + \int_{\partial \Omega} \hat{\varphi}_h  \tilde{\mathbf{v}_h} \cdot \mathbf{n} d x = 0 ,
\end{equation*}
where $\mathbf{n}$ is the edge normal vector and the tilde symbol represents the approximated solution of the corresponding term at the interfaces of elements.

\begin{figure}[ht]
  \centering
    \subfloat[Spherical moving frames]{\label{SphericalMF} \includegraphics[
    width=4.0cm]{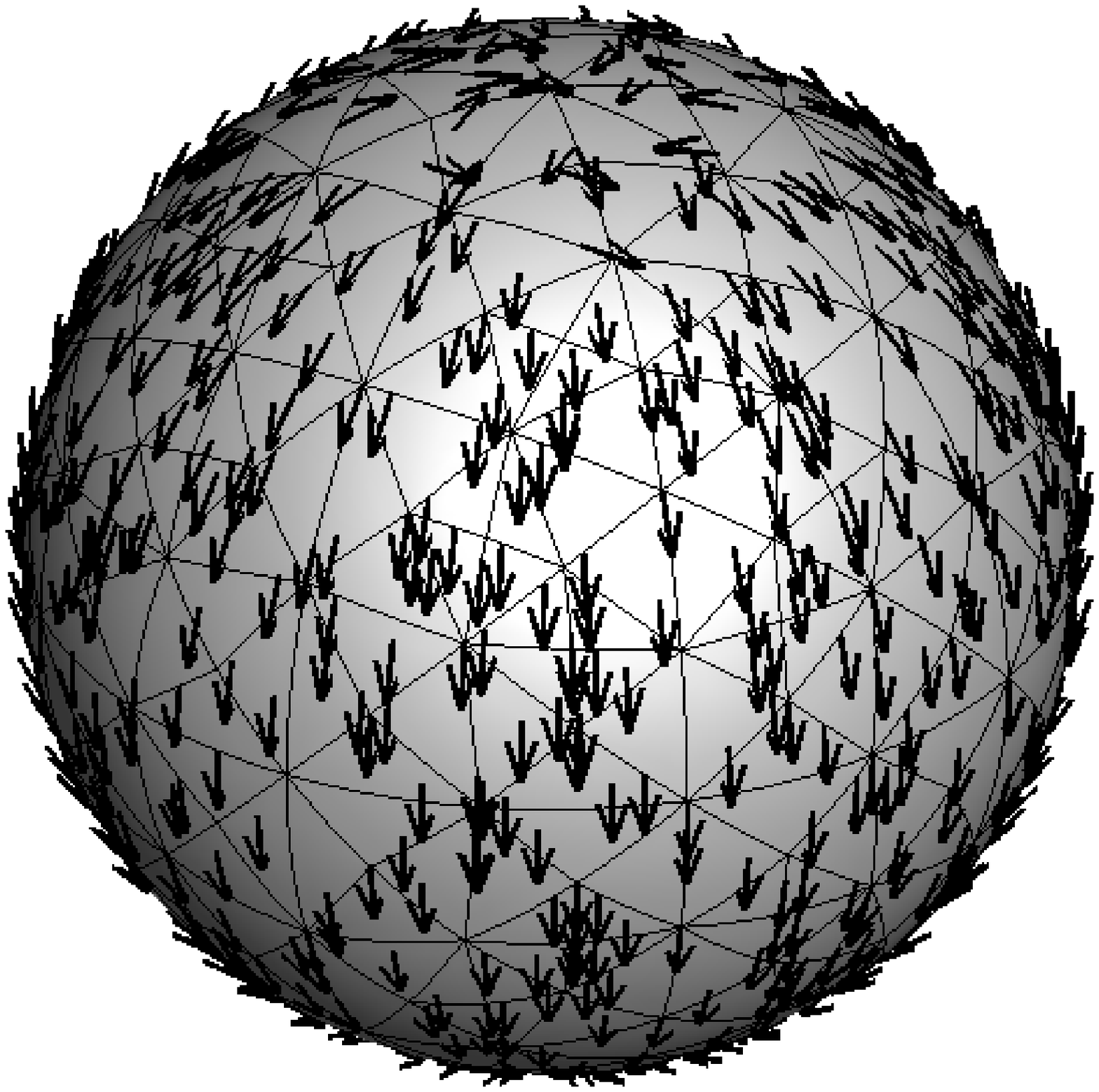}}
  \subfloat[Local moving frames]{\label{LocalMF} \includegraphics[
  width=4.0cm]{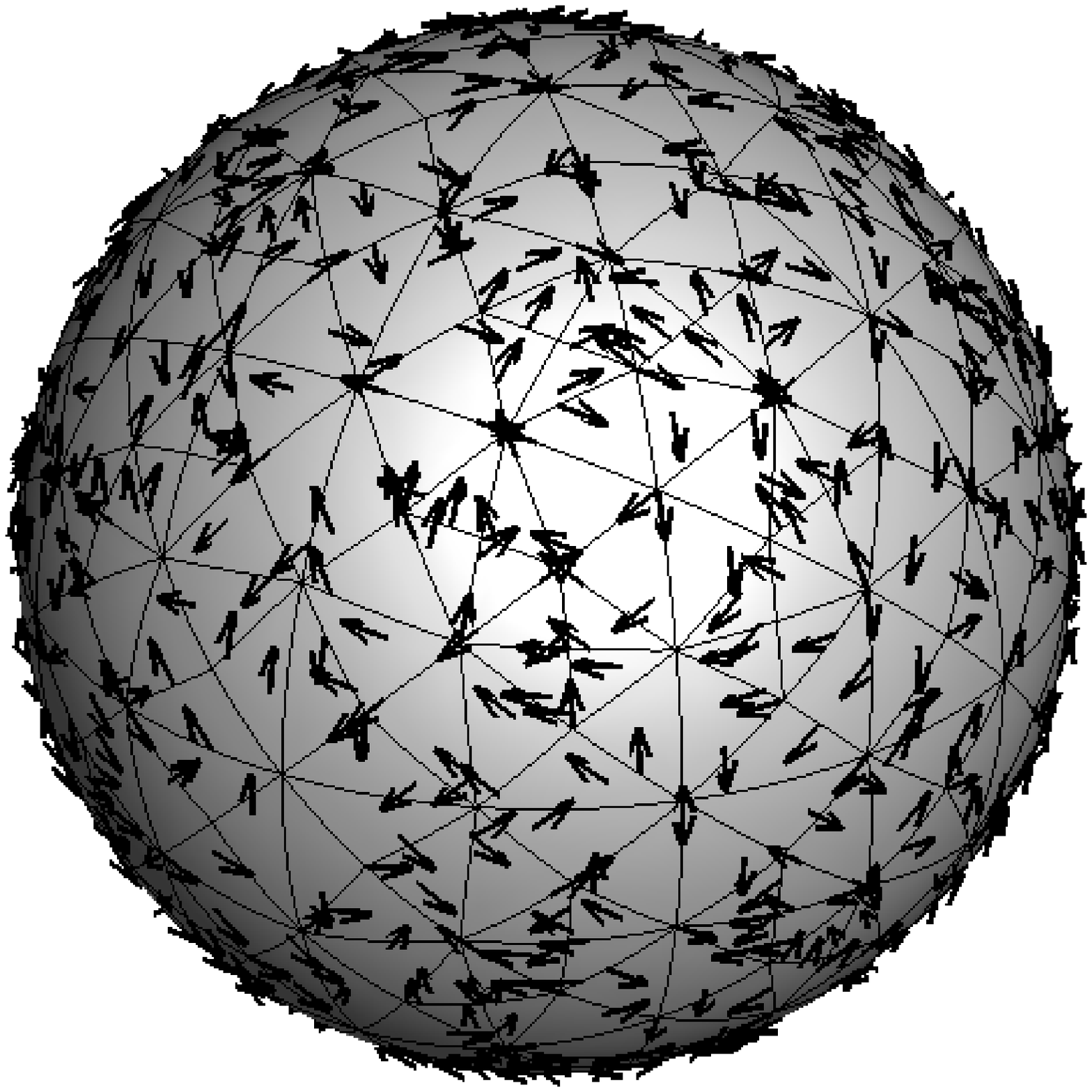}}
  \caption{Distributions of $\mathbf{e}^1$ for Spherical and Local moving frames.}
    \label{MovingFrames}
\end{figure}

\section{Test cases for differential operators}

This section presents the numerical errors of three different types of curved meshes, i.e., ProjMesh, Gmsh, and NekMesh, which are used for the computation of four differential operators, i.e., divergence, gradient, and curl operators. The corresponding test problems are not valid on the entire sphere with the spherical coordinate axis due to the singularities at the poles. Thus, only curved elements, except the small number of elements close to both poles, have the spherical alignment and are used for error computation as shown in Fig. \ref{SphericalMF}. If moving frames are aligned \textit{locally} in each element and possibly discontinuously across the interface of elements, the singularity of the sphere disappears and the entire sphere can be used for error computation, as shown in Fig. \ref{LocalMF}.

For the test problems in this section, the following velocity flow of Rossby-Haurwitz wave is considered: For the spherical shell of the metric $ds^2 = r^2 d \theta^2 + r^2 \sin^2 \theta d \phi^2 $, the velocity flow is defined as \cite{Williamson1992}
\begin{align*}
\mathbf{v} &= v_{\phi}  \hat{\boldsymbol{\phi}} + v_{\theta}  \hat{\boldsymbol{\theta}} ,\\
\mbox{where}&, \\
v_{\phi} &=  \omega \sin \theta + K \sin^3 \theta ( 4 \cos^2 \theta - \sin^2 \theta ) \cos 4 \phi, \\
v_{\theta} &=  -4 K  \sin^3 \theta \cos \theta \sin 4 \phi , 
\end{align*}
where $\omega = K = 7.848 \times 10^{-6}~ s^{-1}$. The four differential operators are evaluated in both direct and weak formulation. The direct method involves the covariant differentiation of the physical space on grid points, and the weak formulation involves performing differentiation, followed by a transformation into modal space. For the full descriptions of direct and weak formulation of each differential operator, refer \cite{MMFCovariant} and \cite{MMF1}, \cite{MMF2}, \cite{MMF3}, \cite{MMF4}, respectively.

\subsection{Divergence}

Let us express the velocity vector on the basis of moving frames such as $\mathbf{v} = v_1 \mathbf{e}^1 + v_2 \mathbf{e}^2$. The divergence of the velocity vector $\mathbf{v}$ on the sphere of radius $r$ is then evaluated as follows \cite{MMFCovariant}.
\begin{equation}
\nabla \cdot \mathbf{v} = \left ( \nabla v_1 \cdot  \mathbf{e}^1 - \Gamma^2_{11} v_2 + \nabla v_2 \cdot  \mathbf{e}^2 + \Gamma^2_{21} v_1  \right ) , \label{DivDirect}
\end{equation}
where $\Gamma^i_{jk}$ is the Christoffel symbol of the second kind and is computed using the connection of $\omega^i_j$ because of the equality $\Gamma^i_{jk} = \omega^i_j \langle \mathbf{e}^k \rangle$ for orthonormal basis $\mathbf{e}^i$. The divergence can be also obtained by a weak formulation with a differentiable test function $\varphi$ in moving frames as follows \cite{MMF1}.
\begin{equation}
\int_{\Omega}  \nabla \cdot \mathbf{v}  \varphi d x = \sum_{m=1}^2 \left [ - \int_{\Omega}  ( \nabla \varphi \cdot  v_m \mathbf{e}^m ) d x + \int_{\partial \Omega}   v_m \mathbf{e}^m \cdot \mathbf{n}  \varphi d s \right ] ,  \label{DivWeak}
\end{equation}
where $\mathbf{n}$ is the edge normal vector. By integration by parts, the covariant derivative is changed to Euclidean inner product and the metric tensor vanishes. The obtained value of divergence is compared to the exact value of the divergence of $\mathbf{v}$ for Rossby-Haurwitz wave which is zero. 


Fig. \ref{DivDirect} presents the convergence of the divergence error for the three different meshes by the direct covariant computation. On the spherical moving frames, all the three meshes show the exponential convergence. Gmsh shows the best accuracy, and Nekmesh shows the least accuracy. This indicates that the differentiation of ProjMesh or Gmsh along the spherical coordinate axis is particularly more accurate than that of NekMesh because corresponding grid points are optimally placed for differentiation along the spherical axis.

For local moving frames, Nekmesh shows the best accuracy, whereas Projmesh and Gmsh error stagnate around $10^{-3}$ when $p \ge 4$. This implies that the grid points Projmesh and Gmsh which are not aligned along the spherical coordinate axis, particularly for inner grid points, are not optimally placed. In other words, geometric approximation error for Projmesh and Gmsh is relatively large along an oblique angle to the spherical coordinate axis. On the other hand, Nekmesh's geometric approximation error for all the direction is trivial enough to show exponential convergence.

In weak formulation, the overall placement of all the grid points within an element is important for the accuracy of derivative regardless of direction, rather than the individual accuracy of each grid point as required for direct formulation. Thus, only Nekmesh for both Spherical and local moving frames shows the exponential convergence without the bottleneck of convergence, as shown in Fig. \ref{DivWeak}. This confirms that Nekmesh's geometric approximation error is trivial not just for grid points along the spherical axis, but also for other grid points inside edges.

\begin{figure}[ht]
  \centering
  \subfloat[Direct formulation]{\label{DivDirect} \includegraphics[
  width=5cm]{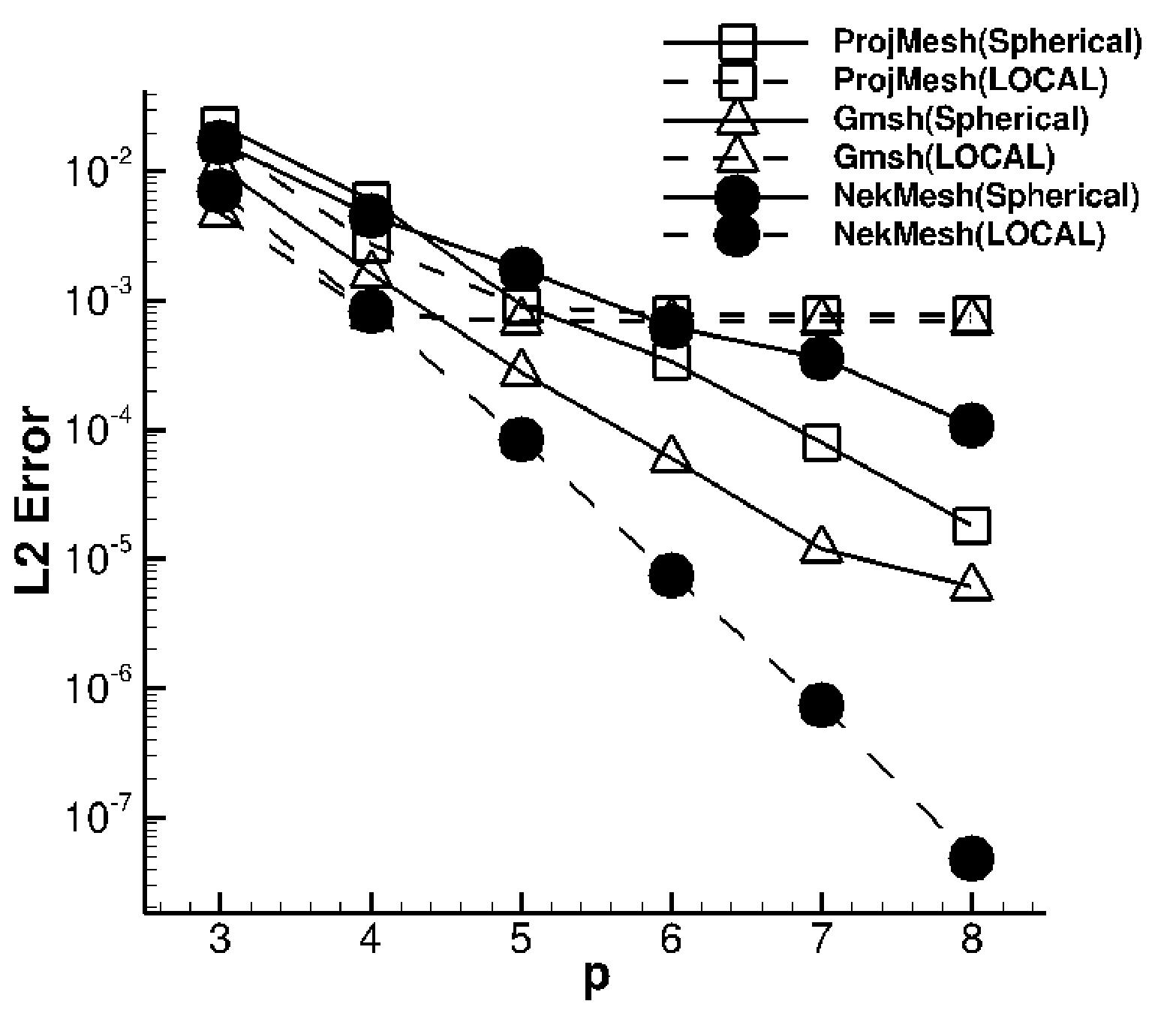}}
  \subfloat[Weak formulation]{\label{DivWeak} \includegraphics[
  width=5cm]{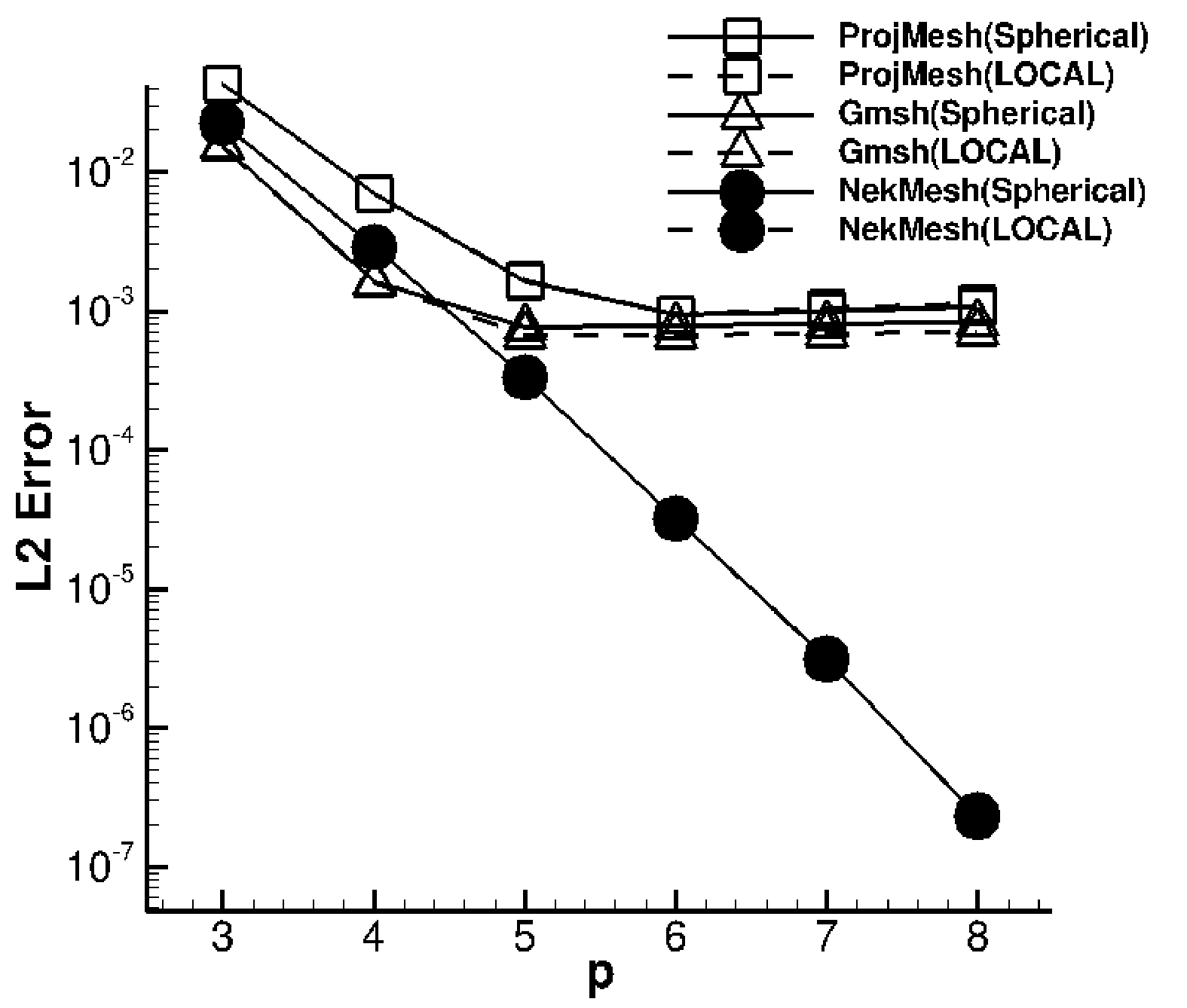}}
    \caption{Divergence error via direct or weak formulation. Solid line = Spherical MF, Dashed line = Local MF, Square = ProjMesh, Triangle = Gmsh, Circle = NekMesh. }
    \label{DivTest} 
\end{figure}

\subsection{Curl}

For a velocity vector $\mathbf{v}$ of the Rossby-Haurwitz flow, the curl of the vector is computed covariantly as follows \cite{MMFCovariant}
\begin{equation}
( \nabla \times \mathbf{v} ) \cdot \mathbf{r}   =   \left ( \nabla v_2 \cdot  \mathbf{e}^1 + \Gamma^2_{12} v_2 + \nabla v_1 \cdot  \mathbf{e}^2 - \Gamma^2_{11} v_1  \right ), \label{CurlDirect} 
\end{equation}
where $\mathbf{r}$ is the radial direction from the center of the sphere. $\Gamma^i_{jk}$ is similarly computed by using the connection 1-form $\omega^i_j$. For a weak formulation for curl operator, the following MMF scheme is used \cite{MMF4}
\begin{align}
\lefteqn{\sum_{i=1}^2 \left [ \int_{\Omega} (\nabla \varphi \cdot \mathbf{e}^{3i} ) v_i d x + \int_{\Omega} v_i \mathbf{e}^i \cdot (\nabla \times \mathbf{e}^3 ) \varphi dx  \right ] } \nonumber \\
& \hspace{4cm} + \int_{\partial {\Omega}} \mathbf{e}^3 \cdot (\mathbf{n} \times \mathbf{v}^* ) \varphi ds = \int_{\Omega} f \varphi d x  , \label{CurlWeak} 
\end{align}
where $\mathbf{v}^*$ is the numerical flux at the interface of the element, and we introduce the new variable $\mathbf{e}^{3i}$ defined as follows $ \mathbf{e}^{3i} = \mathbf{e}^3 \times \mathbf{e}^i$. Note that $\mathbf{e}^3 = \mathbf{r}$. Computed values from Eqs. \eqref{CurlDirect} and \eqref{CurlWeak} are compared to the exact value of the curl of $\mathbf{v}$ such as
\begin{equation*}
(\nabla \times \mathbf{v} ) \cdot \mathbf{r} = - 2 \omega \cos \theta + 30 K \sin^4 \theta \cos \theta \cos 4 \phi .
\end{equation*}

\begin{figure}[ht]
  \centering
  \subfloat[Direct formulation]{\label{CurlConvDirect} \includegraphics[
  width=5.0cm]{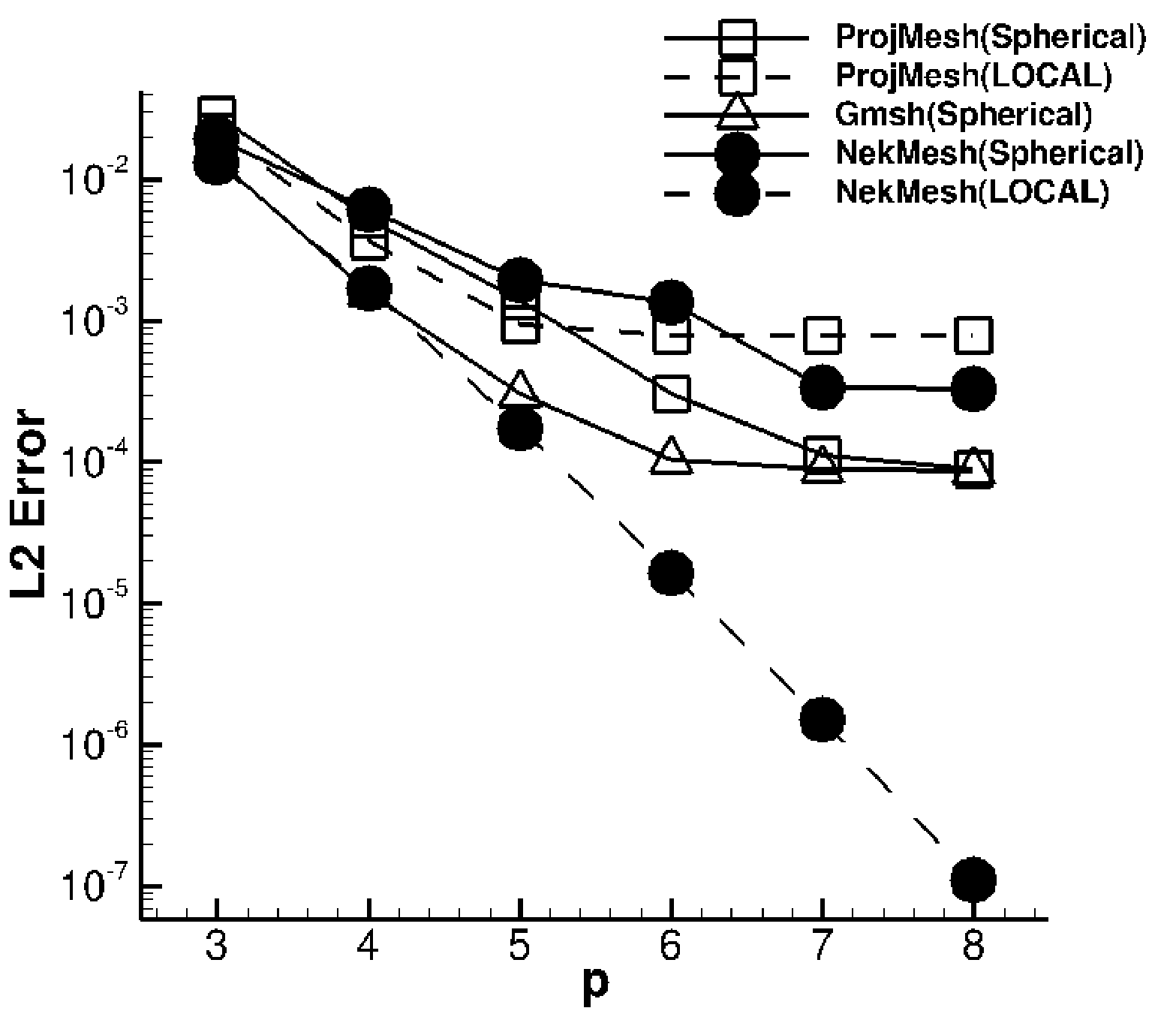}}
  \subfloat[Weak formulation]{\label{CurlConvWeak} \includegraphics[
  width=5.0cm]{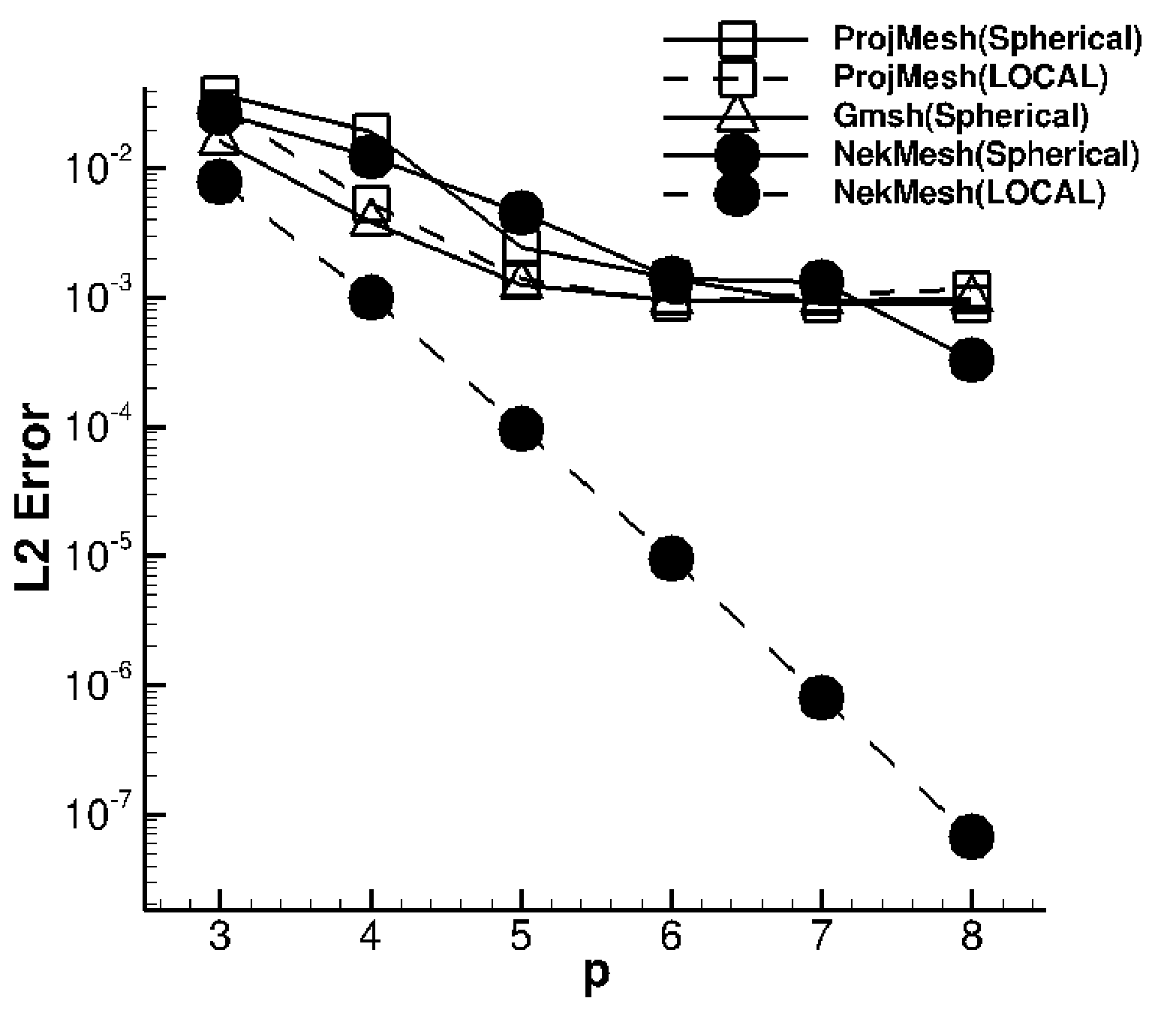}}
    \caption{Curl error via direct or weak formulation. Solid line = Spherical MF, Dashed line = Local MF, Square = ProjMesh, Triangle = Gmsh, Circle = NekMesh. }
    \label{CurlTest}
\end{figure}

Fig. \ref{CurlTest} illustrates similar convergence behavior as that in the case of divergence. Gmsh does not show convergence result for the curl test, thus is omitted in Fig. \ref{CurlTest}. For direct covariant curl computation, all the three meshes shows exponential convergence up to $p=6$, but stagnates from $p=7$. The direct covariant error of ProjMesh with local moving frames stagnates from $p=5$, but the covariant error of Nekmesh shows exponential convergence up to $p=8$. For weak formulation of curl operator, error stagnation is more obvious for ProjMesh and Gmsh for $p \ge 5$. However, Nekmesh shows exponential convergence both for spherical and local moving frames. 

\subsection{Gradient}

\begin{figure}[ht]
  \centering
    \subfloat[Direct formulation]{\label{GradThDirect} \includegraphics[
    width=5cm]{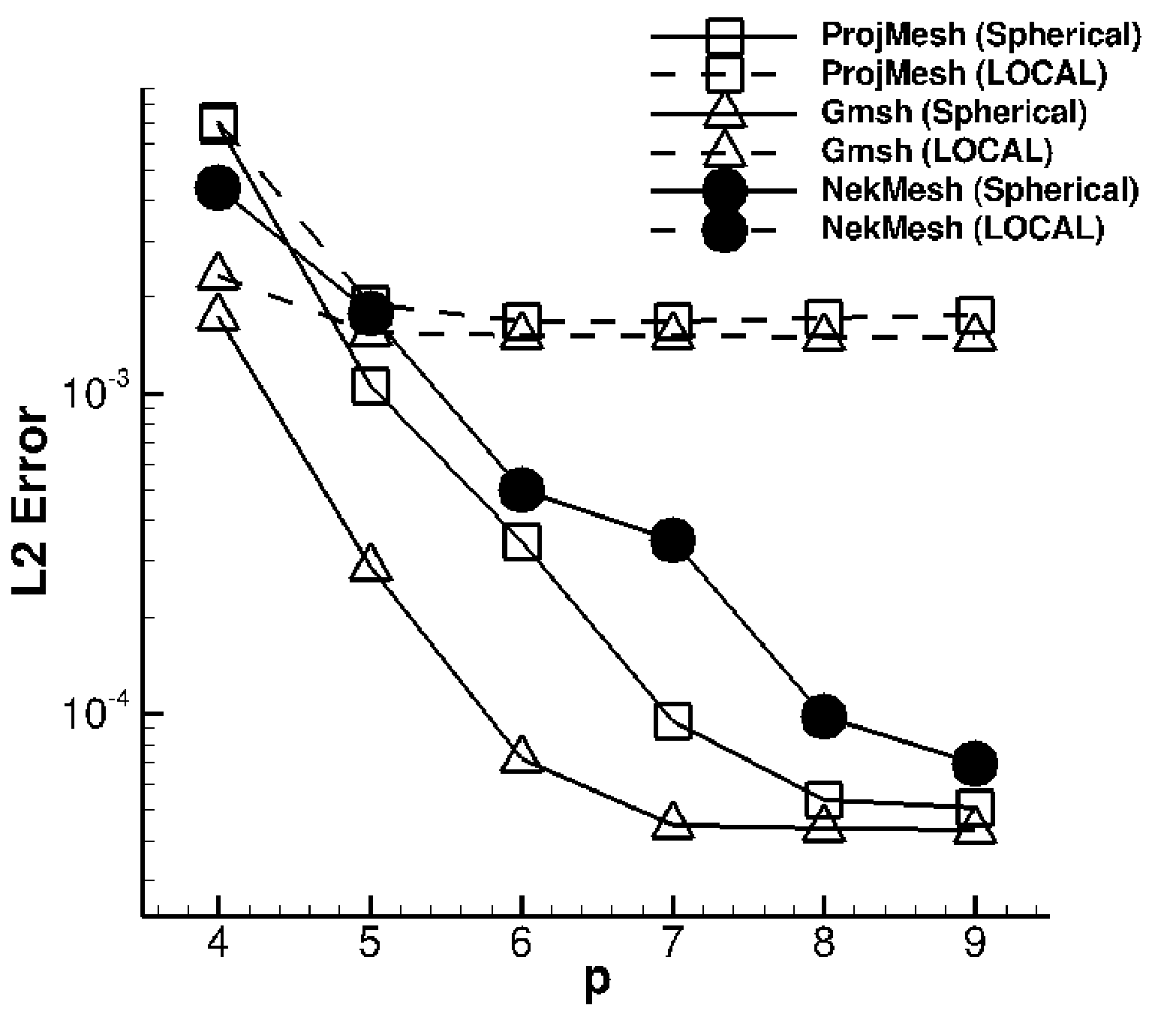}}
  \subfloat[Weak formulation]{\label{GradThWeak} \includegraphics[
  width=5cm]{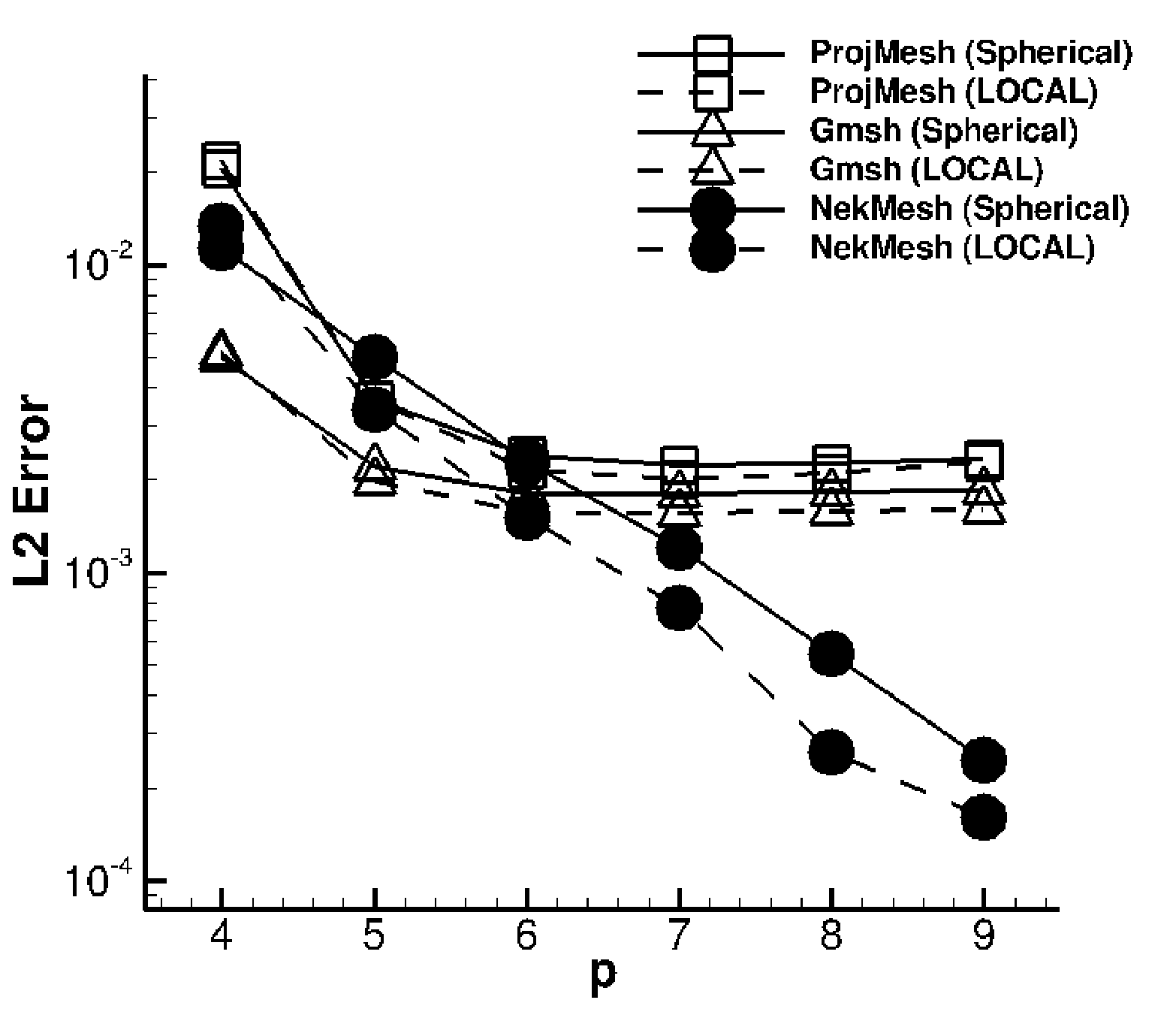}}
  \caption{Convergence of the gradient error via direct or weak formulation. Solid line = Spherical MF, Dashed line = Local MF, Square = ProjMesh, Triangle = Gmsh, Circle = NekMesh.}
    \label{GradTest}
\end{figure}

For a scalar $f$ from the Rossby-Haurwitz wave such as
\begin{equation*}
 f = \omega \sin \theta + K \sin^3 \theta ( 4 \cos^2 \theta - \sin^2 \theta ) \cos 4 \phi,
\end{equation*}
the gradient of $f$ is derived as follows. 
\begin{equation}
\nabla f = \frac{\partial f}{\partial \theta} \boldsymbol{\theta} + \frac{1}{\sin^2 \theta} \frac{\partial f}{\partial \phi} \boldsymbol{\phi}  ,\label{gradf}
\end{equation}
where
\begin{align*}
\frac{\partial f}{\partial \theta} &= \omega \cos \theta + K \sin^2 \theta [ 3 \cos \theta ( 4 \cos^2 \theta - \sin^2 \theta ) - 10 \sin^2 \theta \cos \theta ] \cos 4 \phi, \\
\frac{\partial f}{\partial \phi} &= - 4 K \sin \theta ( 4 \cos^2 \theta - \sin^2 \theta ) \sin 4 \phi.
\end{align*}
By direct covariant differentiation, the gradient can be computed as follows.
\begin{equation}
\nabla f = ( \nabla f \cdot \mathbf{e}^1 ) \mathbf{e}^1 +  ( \nabla f \cdot \mathbf{e}^2 ) \mathbf{e}^2.  \label{CovGrad}
\end{equation}
Note that Eq. \eqref{CovGrad} is the same as the direct covariant formulation of Eq. \eqref{gradf} because the magnitude of the tangent vector of the axis $\boldsymbol{\phi}$ is $\sin \theta$. For weak formulation, the following scheme with moving frames is adapated as follows \cite{MMF3}. For the weak formulation of the gradient of $f$, we have
\begin{equation*}
\int_{\Omega} \nabla f d x = \left ( \int_{\Omega}  \nabla f \cdot \mathbf{e}^1 \varphi dx  \right )  \mathbf{e}^1 +  \left ( \int_{\Omega}  \nabla f \cdot \mathbf{e}^2 \varphi dx  \right )  \mathbf{e}^2,
\end{equation*}
where each component for $1 \le i \le 2$ is computed as follows. 
\begin{equation}
\int_{\Omega} \nabla f \cdot \mathbf{e}^i \varphi dx =   - \left [ \int_{\Omega}  ( \nabla \varphi \cdot \mathbf{e}^i ) f d x + \int_{\Omega} f ( \nabla \cdot \mathbf{e}^i ) \varphi d x \right ] + \int_{\partial \Omega}  ( \mathbf{e}^i \cdot \mathbf{n} ) \varphi \tilde{f} d s ,   \label{GradWeak}
\end{equation}
where $\tilde{f}$ is the numerical flux and is chosen as the upwind flux for $f \mathbf{e}^i$. Fig. \ref{GradTest} confirms the similar results as shown in divergence and curl. In principle, the computation of the gradient should be independent of the direction of differentiation. However, Fig. \ref{GradThDirect} shows the two different convergence for ProjMesh and Gmsh because of nontrivial geometric approximation error other than the spherical coordinate direction. On the other hand, NekMesh shows the same error and convergence rate both for Spherical and Local frames. For weak formulation, the overall geometric approximation error affects the differentiation for any direction, the error of ProjMesh and Gmsh does not converge after $p=5$ as shown in Fig. \ref{GradThWeak}. However, Nekmesh shows exponential convergence both for spherical and local frames. 

\section{Test cases for time-dependent PDEs}

This section presents improvements in the accuracy and conservation properties in the numerical solution of time-dependent partial differential equations on high-order curvilinear mesh, compared with those on low-order curvilinear mesh. To achieve these goals, four PDEs are numerically solved on the sphere: conservation laws, diffusion equations, shallow water equations, and Maxwell's equations. In prior studies \cite{MMF1,MMF2,MMF3,MMF4}, all the meshes of a sphere were generated by the ProjMesh, and therefore, the performance of Nekmesh can only be compared to that of ProjMesh, particularly to understand the correlation between geometric approximation error and conservational properties. The performance of Gmsh is similar to that of ProjMesh, which is the same as the tests for differential operators in previous sections. However, in some tests Gmsh exhibits less stable or less accurate performance compared to ProjMesh. Therefore, the test results of Gmsh are redundant to the conclusion of this section and are not included in the results.

Similar to differential operators, the numerical schemes of each PDE uses moving frames on the surface. Details of each numerical scheme can be found in the corresponding refs. \cite{MMF1,MMF2,MMF3,MMF4}, but a short description will given in each subsections.

\begin{figure}[ht]
  \centering
\subfloat[Overall error]{\label{ErrCase1} \includegraphics[
width=5cm]{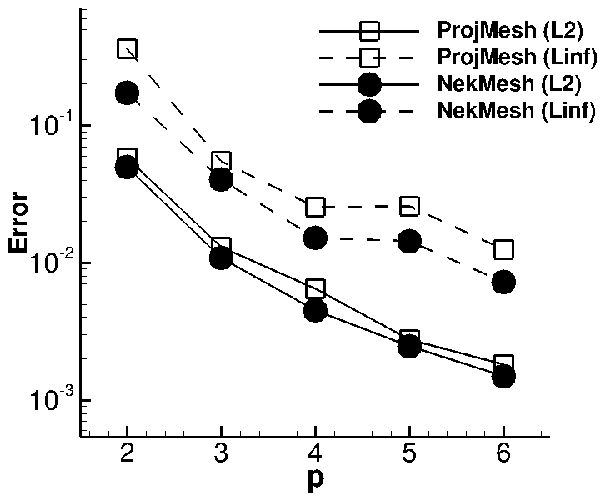}}
   \subfloat[Mass error]{\label{MErrCase1} \includegraphics[
   width=5cm]{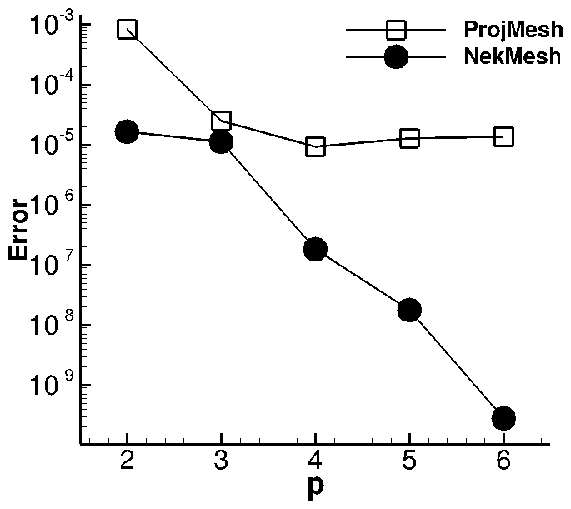}} \\
    \caption{Test problem for conservation laws. $\theta_c=\pi/4$ and a rotational angle of $\pi/2$. $h=0.4$.}
 \label{AdvTest} 
\end{figure}

\begin{table}[ht]
\resizebox{\textwidth}{!}{
\begin{tabular}{c  c c c c c c c c}
 \hline\noalign{\smallskip} 
 $p$ &  \multicolumn{2}{c}{dof} &  \multicolumn{2}{c}{$L_2$ error} & \multicolumn{2}{c}{$L_{\infty}$ error}  &  \multicolumn{2}{c}{Mass error} \\
  \cmidrule(lr){2-3}  \cmidrule(lr){4-5} \cmidrule(lr){6-7} \cmidrule(lr){8-9}
 &  ProjMesh & NekMesh & ProjMesh & NekMesh & ProjMesh & NekMesh & ProjMesh & NekMesh \\
  \noalign{\smallskip}\hline\noalign{\smallskip}
 2 & 5760 & 5976 & 5.84E-02 & 4.97E-02 & 3.62E-01 & 1.73E-01 & 8.27E-04 & 1.66E-05  \\
 3 & 9600 & 9960 & 1.30E-02 & 1.09E-02 & 5.50E-02 & 4.03E-02 & 2.51E-05 & 1.11E-05 \\
 4 & 14400 & 14940 & 6.45E-03 & 4.51E-03 & 2.56E-02 & 1.52E-02 & 9.31E-06 & 1.80E-07 \\
 5 & 20160 & 20916 & 2.78E-03 & 2.48E-03 & 2.57E-02 & 1.43E-02 & 1.26E-05 & 1.75E-08 \\ 
 6 & 26880 & 27888 & 1.83E-03 & 1.50E-03 & 1.24E-02 & 7.27E-03 & 1.34E-05 & 2.82E-10 \\
\end{tabular}}
\caption{Test problem for conservational laws. $\theta_c=\pi/4$ and rotational angle = $\pi/2$. $h=0.4$.}
\label{table::CLtest}
\end{table}

\subsection{Conservation laws}
For a particular distribution of a scalar value $u$ and the velocity vector $\mathbf{v}$, the numerical scheme of the following conservation laws on the sphere
\begin{equation}
\frac{\partial u}{\partial t} + \nabla \cdot (u \mathbf{v} ) =0  \label{CLaws}
\end{equation}
is given as \cite{MMF1}
\begin{equation}
\int\frac{\partial u}{\partial t} \varphi dx =  \sum_{m=1}^2 \left [  \int_{\Omega}  ( \nabla \varphi \cdot   \mathbf{e}^m )  u v^m d x - \int_{\partial \Omega} ( \mathbf{e}^m \cdot \mathbf{n} )  \tilde{u} v^m   \varphi d s \right ] ,  \label{CLscheme}
\end{equation}
where $\tilde{u}$ denotes the upwind flux of the variable $u$. Two moving frames $\mathbf{e}^1$ and $\mathbf{e}^2$ lie on the tangent plane at the point $P$ of the surface. The velocity $\mathbf{v}$ is expanded by the moving frames such as $\mathbf{v} = v^1 \mathbf{e}^1 + v^2 \mathbf{e}^2$. For the numerical test of the scheme on each mesh, the following cosine bell is constructed at the location $(\phi_c, \theta_c)$ in the spherical coordinate system such that
\begin{equation*}
u(\mathbf{x},t=0) = \left \{
\begin{array}{cc}
0.5 ( 1 + \cos ( \pi \mbox{dist}(\varphi,\theta)/r_{\ell}  ))~~~~~~~~~~~~ & \mbox{if}~~ \mbox{dist}(\varphi,\theta) <  r_{\ell} \\
 0.0~~~~~~~~~~~~&  \mbox{if}~~ \mbox{dist}(\varphi,\theta) \ge  r_{\ell} 
\end{array}
\right . \\
\end{equation*}
The cosine bell makes 360 degree rotation along an arbitrary direction and is compared with the original bell \cite{Williamson1992}. The parameter $r_{\ell}$ is the radius of the bell and is set as $7 \pi / 64$. The function $d(\theta,\phi)$ denotes the distance to the center of the bell $(\phi_c,\theta_c)$. The distribution of the cosine bell and the velocity vector that is aligned $\pi / 4$ with respect to the North pole. An explicit fourth-order Runge-Kutta scheme is used for time marching with $\Delta t = 0.0001$ up to the final time $2.0$ with angular frequency of $\pi$.

Fig. \ref{AdvTest} and Table \ref{table::CLtest} illustrate the $L_2$ and $L_{\infty}$ error of the cosine bell after one iteration along the great circle in a specific direction. The test case corresponds to the initial location of the cosine bell at $(\phi_c,\theta_c)=(\pi/4, 3 \pi /4)$ and the velocity vector in the direction of $\pi/4$. The $L_2$ and $L_{\infty}$ errors of Nekmesh are smaller than those of ProjMesh for all $p$. The greatest and smallest difference in $L_2$ is $70~\%$ and $89~\%$, respectively, whereas the greatest and smallest difference in $L_{\infty}$ is $48\%$ and $73\%$, respectively. Note that for this conservation laws test, there is no stagnation of the total error because geometric approximation error prevails, even for ProjMesh. This implies that for this test of conservation laws the geometric approximation error has negligible impact on the overall accuracy.

However, Fig. \ref{MErrCase1} presents that the geometric approximation error has significant impacts on mass conservation. The conservation of the total mass, which is measured as $\int u dx $ on the sphere, should always be preserved because of the conservation properties of Eq. \eqref{CLaws}. However, Fig. \ref{MErrCase1} presents that mass conservation error is not diminishing for $p \ge 3$ by ProjMesh. The mass conservation error stagnates at approximately $10^{-5}$. However, the mass conservation error by NekMesh continues to exponentially converge up to $2 \times 10^{-10}$ at $p=6$.

\begin{figure}[ht]
  \centering
\includegraphics[
width=5cm]{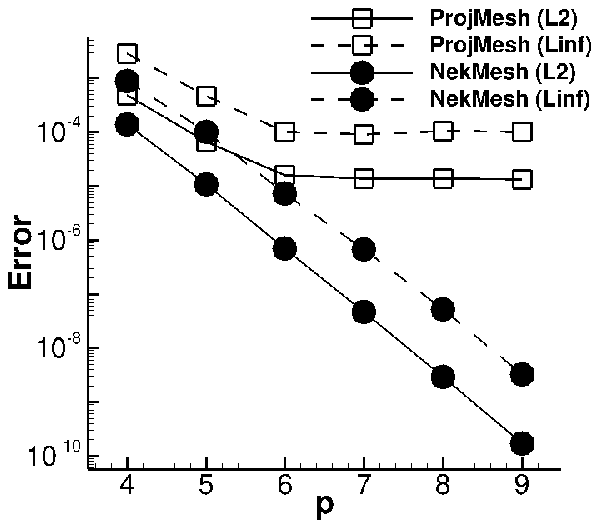}
    \caption{Overall error of the test problem for diffusion-reaction equations on the sphere with $h=0.4$. }
    \label{DiffusionErr} 
\end{figure}

\begin{table}[ht]
\begin{tabular}{c c c c c c c}
 \hline\noalign{\smallskip} 
 $p$  & \multicolumn{2}{c}{dof}  &  \multicolumn{2}{c}{$L_2$ error} & \multicolumn{2}{c}{$L_{\infty}$ error}   \\
 \cmidrule(lr){2-3} \cmidrule(lr){4-5}   \cmidrule(lr){6-7}
 & ProjMesh & NekMesh & ProjMesh & NekMesh & ProjMesh & NekMesh \\
  \noalign{\smallskip}\hline\noalign{\smallskip}
 4 &  14400 & 14940 &  4.77E-04 &    1.37E-04 &  2.82E-03 & 8.97E-04  \\
 5 & 20160 & 20916  &  6.68E-05  &  1.86E-05  & 4.69E-04  &  1.01E-04 \\
 6 & 26880 & 27888   &  1.61E-05  &   7.06E-07  & 1.02E-04 & 7.27E-06   \\
 7 & 34560 & 35856 & 1.40E-05   &  4.62E-08   & 8.98E-05  & 6.77E-07  \\
 8 & 43200 & 44820 & 1.37E-05   &  2.90E-09  & 1.03E-04  &  5.14E-08 \\
 9 & 52800  & 54780  & 1.34E-05   &  1.68E-10  & 9.93E-05  & 3.19E-09 \\
 10 & 63360  & 65736  &  1.32E-05  &  1.87E-11  & 1.09E-04  &  2.82E-10 \\
\end{tabular}
\caption{Test problem for diffusion-reaction equations on the sphere with $h=0.4$.}
\label{table::DiffusionErr}
\end{table}

\subsection{Diffusion Equations}
Consider the following mixed formulation of an the elliptic operator of the scalar variable $u$ with $f \in \mathbb{R}$ and its gradient $\mathbf{q}$ that lies on curved surfaces such that
\begin{equation*}
\nabla \cdot \mathbf{q} ( \mathbf{x} )  = f, ~~~~~\mathbf{q} = \nabla u (\mathbf{x} )  .
\end{equation*}

For isotropic surfaces, the MMF DG scheme of the diffusion equation is given as follows \cite{MMF2}

\begin{align*}
\int_{\Omega} f \varphi dx &= \sum_{m=1}^2 \left [  \int_{\Omega} q_m ( \nabla \varphi \cdot \mathbf{e}^m ) dx + \int_{\partial \Omega} \tilde{q}_m ( \mathbf{n} \cdot \mathbf{e}^m ) \varphi ds  \right ], \\
\int_{\Omega} q_m \varphi dx &= - \int_{\Omega} \nabla \cdot (\varphi \mathbf{e}^m ) + \int_{\partial \Omega} ( \mathbf{n} \cdot \mathbf{e}^m ) \varphi \tilde{u} ds, ~~~~ m = 1,~ 2 ,
\end{align*}
where we expand the gradient $\mathbf{q}$ such that $\mathbf{q} = q_1 \mathbf{e}^1 + q_2 \mathbf{e}^2$. The tilde sign indicates that the corresponding quantity is the numerical flux that is chosen as follows \cite{Castillo2002}
\begin{align*}
\tilde{\mathbf{q}} &= \{\! \{ \mathbf{q} \}\! \} - \boldsymbol{\alpha} [\![u]\!] - {\beta}  [\![\mathbf{q} ]\!], \\
\tilde{u} &= \{\! \{ u \}\! \} + {\beta}  [\![ u ]\!] ,
\end{align*}
where the double braces and a single bracket are defined as $\{\! \{ A \}\! \} \equiv 0.5 (A^+ + A^-)$ and $[\![ A ]\!] \equiv A^- - A^+ $. For the test, $\alpha = 200 \mathbf{n}$ and $\beta = 0.5 $ are used. The superscript (-,+) denotes the field of the corresponding cell and the field of the neighboring cell, respectively. For the test of this elliptic solver, we used the system of linear reaction-diffusion equations for scalar variables $u$ and $v$ such that
\begin{align*}
\frac{\partial u}{\partial t} &= \mu \nabla^2 u + a u + b v,  \\
\frac{\partial v}{\partial t} &=  \nu \nabla^2 v + c u + d v ,
\end{align*}
where the diffusion coefficients are denoted by $\mu=10^{-3}$ and $\nu=2 \times 10^{-3}$. In addition, we let $a=-6$, $b=4$, $c=5$, $d=-4$. The exact solution for the above equations on the sphere is given as follows: for a certain $n$, we have 
\begin{align*}
u (\theta, \varphi, t) &= \sum_{m=-n}^{m=n} e^{\gamma_n t} ( \tilde{A}^m_n \cosh (\delta_n t) + \tilde{B}^m_n \sinh (\delta_n t) ) Y^m_n (\theta, \varphi), \\
v (\theta, \varphi, t) &= \sum_{m=-n}^{m=n} e^{\gamma_n t} ( \tilde{C}^m_n \cosh (\delta_n t) + \tilde{D}^m_n \sinh (\delta_n t) ) Y^m_n (\theta, \varphi).
\end{align*}
The values of each parameter are provided in \cite{MMF2,Pudy2006}. For the time march, the Runge-Kutta 4th order explicit scheme was used with $\Delta t = 0.0001$ up to the final time $1.0$.

Fig. \ref{DiffusionErr} and Table \ref{table::DiffusionErr} present the $L_2$ and $L_{\infty}$ error versus $p$ for the sphere with $h=0.4$. First, for the overall error both in $L_2$ and $L_{\infty}$, Nekmesh exhibits better accuracy. The difference is smaller for lower $p$, but as $p$ becomes higher, the difference becomes larger. This also happens because the stagnation of error occurs for $p \ge 6$ in ProjMesh. The $L_2$ error of Nekmesh is $28~\%$ of that of ProjMesh at $p=4$, but continues to exponentially decrease up to $1.42 \times 10^{-4}~\%$ of the error of ProjMesh at $p$=$8$. Therefore, the geometric approximation error does affect the overall accuracy of $p \ge 6$.

\begin{figure}[ht]
  \centering
     \includegraphics[
     width=6cm]{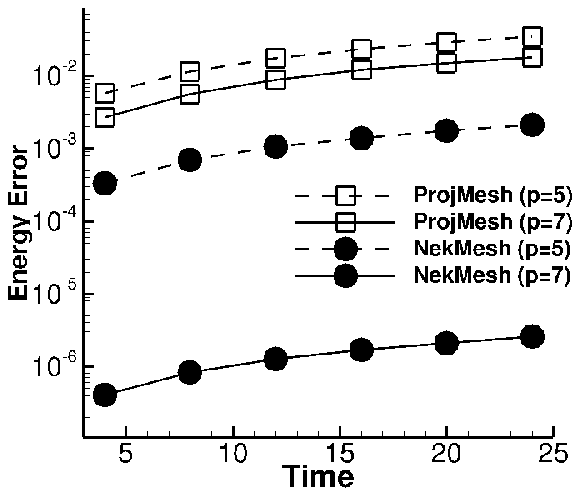} 
    \caption{Energy conservation error of ELF propagation over time up to 24 (18.84 min in real measure).}
         \label{ELFenergy}
    \end{figure}

\begin{table}[ht]
\begin{tabular}{c c c c c}
 \hline\noalign{\smallskip} 
 Time &  \multicolumn{2}{c}{$p=5$} & \multicolumn{2}{c}{$p=7$}   \\
 \cmidrule(lr){2-3} \cmidrule(lr){4-5}  
 & ProjMesh & NekMesh & ProjMesh & NekMesh \\
  \noalign{\smallskip}\hline\noalign{\smallskip}
 4 & 5.76E-03  & 3.36E-04 & 2.69E-03 & 4.05E-07 \\
 8 & 1.14E-02  & 6.94E-04 & 5.55E-03 & 8.24E-07 \\
 12 & 1.74E-02 & 1.06E-03 & 8.87E-03 & 1.28E-06 \\
 16 & 2.34E-02 & 1.39E-03 & 1.19E-02 & 1.67E-06 \\
 20 & 2.90E-02 & 1.74E-03 & 1.48E-02 & 2.10E-06 \\
 24 & 3.44E-02 & 2.11E-03 & 1.78E-02 & 2.57E-06 
\end{tabular}
\caption{Energy conservation error of ELF propagation. $h=0.4$. Dof of ProjMesh for $p$=$5$ and $p$=$7$ is 21060 and 34560. Dof of NekMesh for $p$=$5$ and $p$=$7$ is 20916 and 35858.}
\label{table::ELFenergy}
\end{table}

\subsection{Maxwell's equations}

The third partial differential equation tested on the sphere is the time-dependent Maxwell's equation without source terms such as
\begin{equation*}
\hat{\varepsilon} \frac{\partial \mathbf{E}}{\partial t} = \nabla \times \mathbf{H},  ~~~~\hat{\mu} \frac{\partial \mathbf{H}}{\partial t} = - \nabla \times \mathbf{E} ,  
\end{equation*}
where $\mathbf{E}$ and $\mathbf{H}$ denote the electric field and $H$-field, respectively; whereas $\hat{\varepsilon}$ and $\hat{\mu}$ denote the permittivity and permeability tensor, respectively. The MMF scheme for solving the above system of linear equations was proposed in ref. \cite{MMF4} as follows: Suppose that the domain $\mathcal{M}$ is stationary such that $\partial \mathbf{e}^m / \partial t \cdot \mathbf{e}^i = 0$ for all $i$ and $m$. By introducing the new variables $\mathbf{e}^{ik} \equiv \mathbf{e}^i \times \mathbf{e}^k$, $\mathbf{e}^{i3} = - \mathbf{e}^{3i}$, the weak form of the above equations for the transverse magnetic (TM) mode are obtained for the test function $\varphi$ such that for $i=1,2$, 
\begin{align}
\lefteqn{\int_{\Omega} \mu^i \frac{\partial H^i }{\partial t} \varphi dx  + \int_{\Omega}  \sigma^{*i} H^i  dx -  \int_{\Omega} E^3 \nabla \varphi \cdot  \mathbf{e}^{3i} dx} \hspace{1cm}  \nonumber \\
&    +  \int_{\Omega}  E^3 \mathbf{e}^3 \cdot  ( \nabla \times \mathbf{e}^i ) \varphi dx    + \int_{\partial \Omega} \mathbf{e}^i \cdot \left ( \mathbf{n} \times \mathbf{E}^{3*} \right ) \varphi ds  = 0  ,  \label{TMh}
\end{align}

\begin{align}
\lefteqn{\int_{\Omega} \varepsilon^3 \frac{\partial E^3}{\partial t} \varphi dx + \int_{\Omega}  \sigma^{3} E^3  dx -   \int_{\Omega} H^m \nabla \varphi \cdot  \mathbf{e}^{3m} dx } \hspace{1cm} \nonumber \\
&-  \int_{\Omega}  H^m \mathbf{e}^m \cdot \left ( \nabla \times \mathbf{e}^3 \right ) \varphi d x  - \int_{\partial \Omega} \mathbf{e}^3 \cdot \left ( \mathbf{n} \times  \mathbf{H}^* \right ) \varphi ds = 0 , \label{TMe}
\end{align}
where $\mathbf{n}$ is the edge normal vector. $\mathbf{E}^{3*}$ and $\mathbf{H}^*$ denote the upwind flux at the interface of curved elements defined as
\begin{align}
- \mathbf{e}^i \cdot (\mathbf{n} \times \mathbf{E}^{3*} ) &= \mathbf{e}^i \cdot   \frac{  \left ( - \mathbf{n} \times \mathbf{e}^3 \right ) \{ \! \{ Y_i E^3 \}\! \}+ 0.5 \alpha  \mathbf{n} \times ( \mathbf{n} \times [\![ \mathbf{H} ]\!]  )  }{ \{ \! \{ Y_i \}\! \} }  , \label{TMflux1} \\
 \mathbf{e}^3 \cdot (\mathbf{n} \times \mathbf{H}^* ) &=  \mathbf{e}^3 \cdot  \frac{ n^m \mathbf{e}^m \times \{\!  \{ Z_m \mathbf{H} \}\!  \} - 0.5 \alpha [\![ E^3 ]\!]  }{ \{\!  \{ Z_m \}\!  \} }   ,  \label{TMflux2}
\end{align}
where we introduced the new variable $Z^{\pm}_i = \sqrt{\mu^{3-i} / \varepsilon^3 } = ( Y^{\pm}_i )^{-1}$ such that $\hat{\mu} \mathbf{H} = \mu^1 H^1 \mathbf{e}^1 + \mu^2 H^2 \mathbf{e}^2$. The parameter $\alpha$ is in the range of $0 < \alpha \le 1$. An explicit fourth-order Runge-Kutta scheme is used for time marching with $\Delta t = 0.001$ up to the final time $24.0$.

For the test of Maxwell's equations on the sphere, an extremely low frequency wave propagation is simulated \cite{Simpson2002}. A Gaussian radius of 0.2 (approximately 1270 km in real measure) of electromagnetic impulse is initiated at a point on the sphere and it propagates throughout the surface to reach the antipode of the initiating point and travels backward to propagate back and forth between the initial point and the antipode. The electromagnetic energy propagates along the electromagnetic wave but must be preserved over the time because it is assumed that electromagnetic fields are not dissipative.

The first validation of this numerical test was to confirm the field distribution at the antipode. The first electric field is the positive field followed by the negative field. However, the following electric field is in the reverse order, the negative field is followed by the positive field. The quantification of the test is obtained by computing the total electromagnetic energy conservation loss. Fig. \ref{ELFenergy} and Table \ref{table::ELFenergy} present the energy loss error up to $T$=$12$ (9.42 min in real measure) for ProjMesh and Nekmesh when $h$=$0.4$ with $p$=$5$ or $p$=$7$. For both $p$, it is again confirmed that the electromagnetic energy loss of Nekmesh is significantly smaller than that of Projmesh at all times. We observe that, for $p$=$5$ and $p$=$7$, the energy conservation error of Nekmesh is $11.0~\%$ and $0.41~\%$ times smaller than that of ProjMesh, respectively.

\subsection{Shallow Water Equations}

The fourth partial differential equations tested on the sphere is the shallow water equations given as follows.
\begin{align}
& \frac{\partial H}{\partial t} + \nabla \cdot ( H \mathbf{u}  ) = 0 \,,  \label{SWE1}\\
& \frac{\partial \left( H \mathbf{u} \right)}{\partial t} + \nabla \cdot (H \mathbf{u} \mathbf{u} ) + \frac{g}{2} \nabla H^2 =  f H \left (  \mathbf{u}  \times \mathbf{k} \right ) + g H \nabla H_0\,, \label{SWE2}
\end{align}
where $\eta$ denotes the free surface elevation, $\mathbf{u}$ denotes the depth-averaged velocity, $H_0$ denotes the still water depth, and $H= H_0 + \eta$ denotes the total water depth. $f$ and $g$ correspond to the Coriolis parameter and the gravitational constant, respectively. By expressing $\mathbf{u} = u_1 \mathbf{e}^1 + u_2 \mathbf{e}^2$ on the sphere, the MMF-SWE scheme on an arbitrarily rotating curved surface is derived as follows \cite{MMF3}: For $i=1,2$,

\begin{align}
&\int_{\Omega} \frac{\partial H}{\partial t} \varphi dx  - \int_{\Omega}  H  \mathbf{u} \cdot  \nabla \varphi  dx + \int_{\partial \Omega} \widetilde{H} \tilde{\mathbf{u}} \cdot {\mathbf{n}} \varphi ds = 0 \,, \label{SWEint1} \\
 & \int_{\Omega} \frac{\partial H u_i  }{\partial t} \varphi dx + \int_{\Omega}   H  \left ( u_1  \frac{\partial \mathbf{e}^1}{\partial t} + u_2  \frac{\partial \mathbf{e}^2}{\partial t}  \right ) \cdot \mathbf{e}^i \varphi \,dx  \nonumber \\
 & -  \int_{\Omega}  \left ( H {u}_i^2 +   \frac{g H^2}{2}    \right ) \nabla \varphi \cdot \mathbf{e}^i d x   +  \int_{\Omega}  \left [ \widetilde{H}  \widetilde{u}_i \mathbf{e}^i \cdot \mathbf{n} + \frac{g}{2} \widetilde{H}^2  \right ] \varphi dx  \nonumber \\ 
&  - \int_{\Omega}   H  u_1  u_2  \nabla \varphi \cdot \mathbf{e}^{3-i} dx  - \int_{\Omega}   \frac{g H^2}{2 }  \left ( \nabla \cdot \mathbf{e}^i \right )  \varphi \,dx   \nonumber   \\
   &   = (-1)^{i+1} \int_{\Omega}  f \left( H u_{3-j} \right) \varphi dx + \int_{\Omega}  g H \nabla H_0 \cdot  \mathbf{e}^i \varphi\, dx   .\label{SWEint2} 
 \end{align}

For the above SWE equations, the classical five Williamson's tests are tested on ProjMesh and NekMesh. The Lax-Friedrich flux is used for numerical flux across the interfaces of elements and an explicit fourth-order Runge-Kutta scheme is used for time marching.

\subsubsection{Steady zonal flow}

The first SWE test on the sphere focuses on the steady zonal flow that does not change the free surface elevation ($\eta$) and the velocity vector ($\mathbf{u})$ over time \cite{Williamson1992}. For the horizontal velocity vector $\mathbf{u} = (u_x, u_y, u_z)$ in dimensionless form decomposed into $\tilde{\mathbf{u}} = {u}_{\phi}   \boldsymbol{\phi} +  {u}_{\theta}  \boldsymbol{\theta}$, each component of the velocity vector is given as follows.
\begin{align*}
{u}_{\phi}  &= {u}_0 ( \cos \theta \cos \alpha + \sin \theta \cos \phi \sin \alpha ),  \\
{u}_{\theta}  &= - {u}_0 \sin \phi \sin \alpha,   
\end{align*}
where the magnitude of the initial velocity vector is ${u}_0 = 2 \pi / 12$ and $\alpha$ is the angle of the velocity vector with respect to the polar axis. The free surface elevation ${\eta}$ is given by
\begin{equation*}
\eta = H_0 - \frac{1}{\tilde{g}} \left ( \Omega {u}_0 + \frac{ u^2_0}{2} \right ) ( - \cos \phi \cos \theta \sin \alpha + \sin \theta \cos \alpha )^2 ,   
\end{equation*}
where $\tilde{g}$ is the gravitational constant in dimensionless form derived as $\tilde{g} = g / r_a $ with the radius of the earth ($r_a$) being $6.37122 \times 10^6~ (m)$ and the gravitational constant ($g$) being 9.80616 $m/s^2$. The still water depth $H_0$ is $g H_0 = 2.94 \times 10^4 ~(m^2 /s^2)$ and the angular frequency of the rotation of earth $\Omega$ is $7.292 \times 10^{-5} sec^{-1}$. In addition, the Coriolis parameter is provided such as
\begin{equation*}
{f} = 2 {\Omega} ( - \cos \phi \cos \theta \sin \alpha + \sin \theta \cos \alpha ) .  
\end{equation*}
Substituting these values for SWE equations \eqref{SWE1} - \eqref{SWE2}, or in weak formulations \eqref{SWEint1} - \eqref{SWEint2}, the right hand side of the equations should be zero, which implies that $H$ or $\mathbf{u}$ should remain the same at all time. For the time marching, we set $\Delta t = 0.0005$.

\begin{figure}[ht]
  \centering
  \subfloat[Overall error]{\label{ZonalL2} \includegraphics[
  width=5cm]{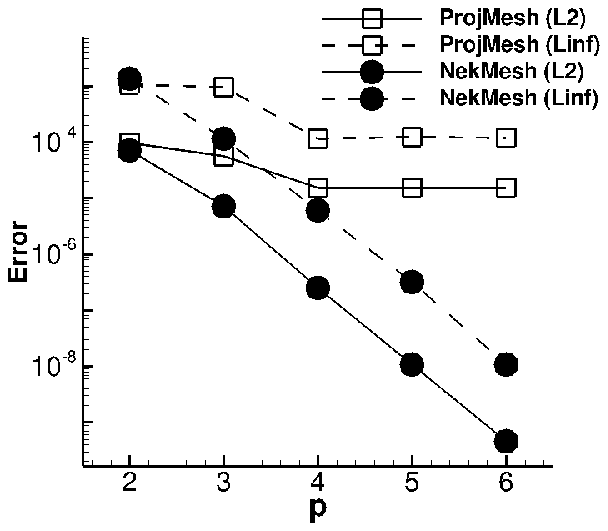}}
  \subfloat[Mass $\&$ energy error]{\label{ZonalMass} \includegraphics[
  width=5cm]{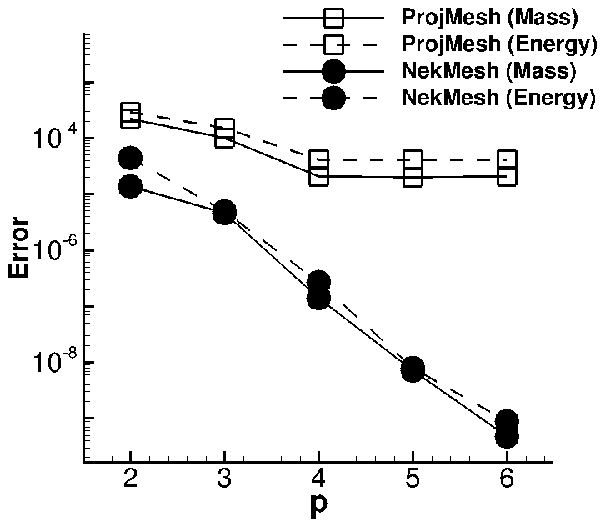}}
    \caption{Steady zonal flow at 5 days. $h$=$0.4$. Error is measured in normalized variables.}
    \label{SWEZonalTest1} 
\end{figure}

Fig. \ref{ZonalL2} and Table \ref{table::SteadyETA} presents the $L_2$ error (solid line) and $L_{\infty}$ error (dashed line) at 5 days for time steps of $0.0005$ (43.2 sec in real measure) and $\alpha = \pi/4$. For ProjMesh, the error does not converge exponentially for $p\ge 4$ and the overall error stagnates by approximately $1.5 \times 10^{-5}$ in $L_2$ and $1.1 \times 10^{-4}$ in $L_{\infty}$. Moreover, the convergence rate for ProjMesh is not at an exponential rate even for lower $p$. From this convergence graph, we conjecture that the geometric approximation error of ProjMesh at the order of $ 10^{-7}$ contributes significantly to the overall $L_2$ and $L_{\infty}$ error after the long time integration up to $T=5.0$. On the other hand, the error in Nekmesh continues to converge up to $p=6$: $4.6 \times 10^{-10}$ in $L_2$ and $1.1 \times 10^{-8}$ in $L_{\infty}$.

In addition, geometric approximation error seems to significantly affect the energy and mass conservation properties, similar to conservation laws. Fig. \ref{ZonalMass} and Table \ref{table::SteadyETA} present the relative loss error of mass and energy to the total mass and energy. This figure confirms that Nekmesh exhibits superior performance in conservation properties such as mass and energy as in shallow water equations. The mass and energy error of Nekmesh converges exponentially up to $10^{-9}$, but those of ProjMesh fail to converge for $p \ge 4$ and stagnates at approximately $2 \times 10^{-5}$.

\begin{figure}[ht]
  \centering
\  \subfloat[$p=2$]{\label{ZonalComp1} \includegraphics[
width=5cm]{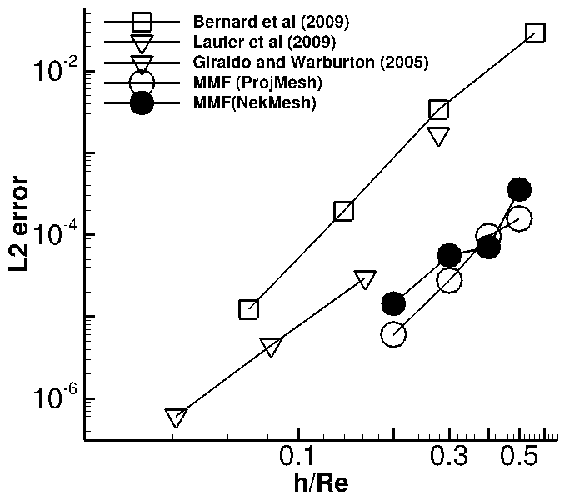}}
  \subfloat[$p=4$]{\label{ZonalComp2} \includegraphics[
  width=5cm]{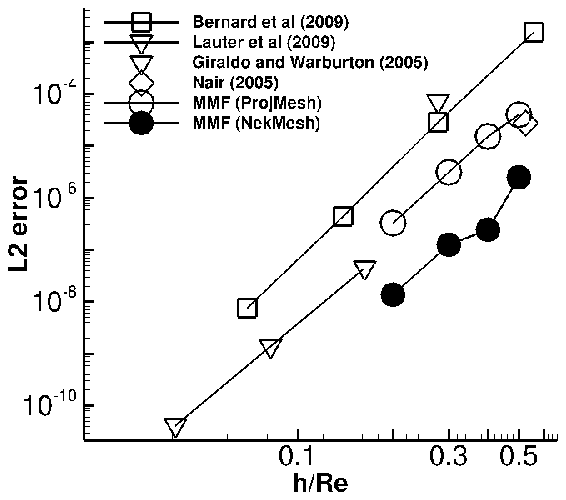}}
    \caption{Comparisons of the normalized $L_2$ error for water height versus $h/R$. }
    \label{SWEZonalTComp} 
\end{figure}

Fig. \ref{SWEZonalTComp} demonstrates the competitive performance of the MMF-SWE scheme with ProjMesh and NekMesh. When a lower order of $p=2$ is used, the accuracy of ProjMesh and Nekmesh is almost indistinguishable. However, for a higher p, i.e. $p$=$4$, Nekmesh exhibits the superior accuracy, approximately $1 \% \sim 6 \%$ of the overall $L_2$ error of ProjMesh for $ 0.2 \le h \le 0.5$. The MMF-SWE scheme then exhibits the best accuracy among the previously proposed schemes for the steady-state zonal flow test \cite{Legat,Lauter2008,Giraldo2005,Nair2005} . 

\begin{table}[ht]
\resizebox{\textwidth}{!}{
\begin{tabular}{c c c c c c c c c}
 \hline\noalign{\smallskip} 
$p$ &  \multicolumn{2}{c}{$L_2$ error} & \multicolumn{2}{c}{$L_{\infty}$ error}  &  \multicolumn{2}{c}{Mass error}  &  \multicolumn{2}{c}{Energy error} \\
 \cmidrule(lr){2-3} \cmidrule(lr){4-5} \cmidrule(lr){6-7} \cmidrule(lr){8-9}
 & ProjMesh & NekMesh & ProjMesh & NekMesh & ProjMesh & NekMesh & ProjMesh & NekMesh \\
  \noalign{\smallskip}\hline\noalign{\smallskip}
 2 & 9.57E-05 & 7.04E-05 & 1.05E-03 & 1.34E-03 & 2.16E-04 & 1.38E-05 & 2.85E-04 & 4.44E-05 \\
 3 & 5.65E-05 & 7.11E-06 & 9.42E-04 & 1.16E-04 & 9.81E-05 & 4.58E-06 & 1.52E-04 & 4.89E-06 \\
 4 & 1.51E-05 & 2.48E-07 & 1.15E-04 & 5.94E-06 & 2.05E-05 & 1.39E-07 & 4.12E-05 & 2.67E-07 \\
 5 & 1.52E-05 & 1.09E-08 & 1.22E-04 & 3.15E-07 & 1.99E-05 & 7.19E-09 & 4.06E-05 & 7.93E-09 \\
 6 & 1.52E-05 & 4.63E-10 & 1.17E-04 & 1.06E-08 & 2.06E-05 & 4.79E-10 & 4.14E-05 & 8.73E-10 \\
\end{tabular}}
\caption{Test problem for steady zonal flow at 5 days. $h$=$0.4$.}
\label{table::SteadyETA}
\end{table}

\subsubsection{Unsteady zonal flow}

\begin{figure}[ht]
  \centering
  \subfloat[Overall error]{\label{UnsteadyL2} \includegraphics[
  width=5cm]{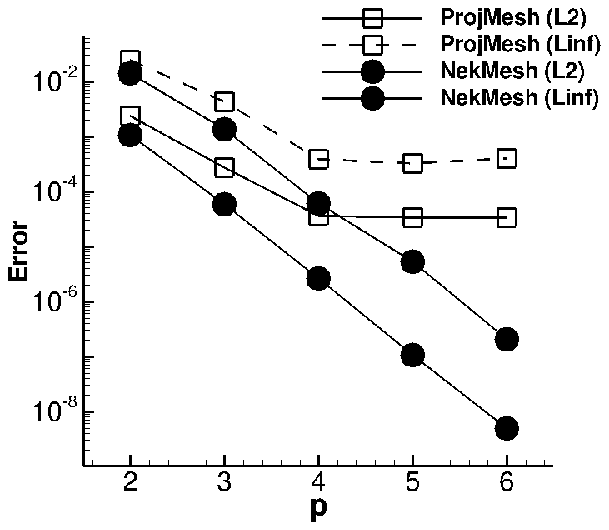}}
  \subfloat[Mass $\&$ energy error]{\label{UnsteadyMass} \includegraphics[
  width=5cm]{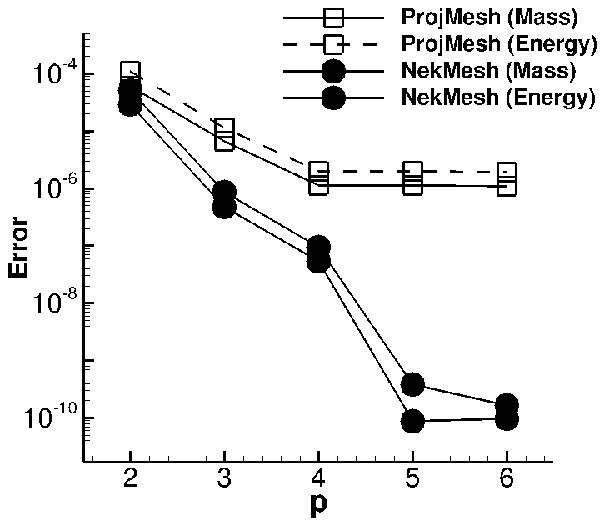}}
    \caption{Unsteady zonal flow at 0.5 days. $h$=$0.4$.}
    \label{SWEUnsteady} 
    \end{figure}

Compared to the steady zonal flow, the unsteady zonal flow has a different magnitude for $\eta$. Thus, the direction of the velocity vector $\mathbf{u}$ depends on the time variable mainly because of non-constant $H_0$ on the sphere. Here are the details of the initial conditions. For the Coriolis parameter $ {f} = 2  {\Omega} \sin \theta$, the velocity vector is given as
\begin{align*}
 {u}_{\phi} ( \phi, \theta, t) &=  {u}_0 ( T_R \sin \alpha \sin \theta + \cos \alpha \cos \theta ) ,  \\
 {u}_{\theta}  ( \phi, t)  &= -  {u}_0  (\sin \phi \cos  {\Omega} t + \cos \phi \sin  {\Omega} t ) \sin \alpha ,
\end{align*}
where  $T_R   =  \cos \phi \cos  {\Omega} t - \sin \phi \sin  {\Omega} t$. We then obtain
\begin{align*}
 {\eta}   ( \phi, \theta, t)  &=   \frac{1}{2  \tilde{g}} \left [ -\left \{   {u}_0 ( - T_R \sin \alpha \cos \theta  + \cos \alpha \sin \theta ) +  {\Omega} \sin \theta   \right  \}^2  + (  {\Omega} \sin \theta )^2 \right ], \\
 {H}_0 (\theta) &=  \frac{133681}{ r_a \tilde{g}} -  \frac{10 }{ r_a \tilde{g}}  - \frac{1}{2  \tilde{g} } (  {\Omega} \sin \theta )^2   .
\end{align*}
For the time marching, we set $\Delta t = 0.0005$. $L_2$ and $L_{\infty}$ errors are computed at $0.5$ days with time steps of $0.0005$ (43 sec in real measurement). Fig. \ref{SWEUnsteady} and Table \ref{table::UnSteadyETA} present the similar accuracy and convergence behavior of NekMesh and ProjMesh compared to the steady-state zonal test; the overall error of NekMesh converges exponentially up to $5.0 \times 10^{-9}$ in $L_2$ and $2.1 \times 10^{-7}$ in $L_{\infty}$. For $p \ge 4$, the overall error of ProjMesh does not decrease below $3.3 \times 10^{-5}$ in $L_2$ and $4.0 \times 10^{-4}$ in $L_{\infty}$. For ProjMesh, the conservation error does not decrease below $1.1 \times 10^{-6}$ for mass and $1.9 \times 10^{-6}$ for energy. However, the conservation error in Nekmesh converges up to $9.7 \times 10^{-11}$ and $1.7 \times 10^{-10}$ for mass and energy, respectively.

\begin{table}[ht]
\resizebox{\textwidth}{!}{
\begin{tabular}{c c c c c c c c c}
 \hline\noalign{\smallskip} 
$p$ &  \multicolumn{2}{c}{$L_2$ error} & \multicolumn{2}{c}{$L_{\infty}$ error}  &  \multicolumn{2}{c}{Mass error}  &  \multicolumn{2}{c}{Energy error} \\
 \cmidrule(lr){2-3} \cmidrule(lr){4-5} \cmidrule(lr){6-7} \cmidrule(lr){8-9}
 & ProjMesh & NekMesh & ProjMesh & NekMesh & ProjMesh & NekMesh & ProjMesh & NekMesh \\
  \noalign{\smallskip}\hline\noalign{\smallskip}
 2 & 2.40E-03  & 1.04E-03  & 2.41E-02  & 1.40E-02  & 6.10E-05  & 2.93E-05  & 1.10E-04 & 5.04E-05  \\
 3 & 2.73E-04  & 5.87E-05  & 4.26E-03  &  1.33E-03 & 6.62E-06  & 4.68E-07  & 1.32E-05 & 8.76E-07  \\
 4 & 3.55E-05  & 2.61E-06  & 3.80E-04  & 6.16E-05  & 1.14E-06  & 5.45E-08  & 2.01E-06 & 9.36E-08  \\
 5 & 3.32E-05  & 1.04E-07  & 3.20E-04  & 5.19E-06  & 1.12E-06  & 8.67E-11  & 2.00E-06 & 3.85E-10  \\
 6 & 3.32E-05  & 5.01E-09  & 3.99E-04  & 2.08E-07  & 1.08E-06  & 9.73E-11  & 1.02E-06 & 1.67E-10  \\
\end{tabular}}
\caption{Test problem for unsteady zonal flow at 0.5 days. $h$=$0.4$ }
\label{table::UnSteadyETA}
\end{table}
   
    \subsubsection{Rossby-Haurwitz flow}
\begin{figure}[ht]
  \centering
      \subfloat[RH wave]{\label{ RossbyFinal } \includegraphics[
      width=5cm]{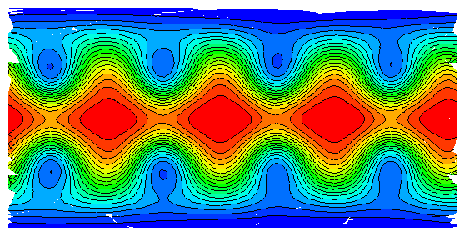}}
  \subfloat[Disturbed RH wave]{\label{DRossbyFinal} \includegraphics[
  width=5cm]{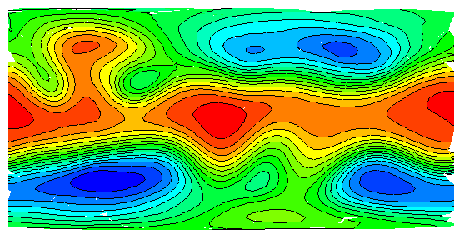}}
    \caption{Rossby--Haurwitz wave and disturbed Rossby--Haurwitz wave at 15th days with $h=0.4$ and $p=5$.}
    \label{SWERossbyFinal} 
    \end{figure}

  The third SWE test problem that appears in the Williamson's classical test problem \cite{Williamson1992} is the Rossby--Haurwitz flow on the sphere. A periodic distribution of $\eta$ and $\mathbf{u}$ in the direction of longitudinal axis slowly rotates around the axis of the sphere. Because there is no exact solution to this test problem, the mass and energy conservation are considered to quantify the quality of the numerical scheme. Consider the velocity components $\mathbf{u} = u \boldsymbol{\phi} + v \boldsymbol{\theta}$ as
\begin{align}
u &=   {\omega} \cos \theta +  {K} \cos^{R-1} \theta ( R \sin^2 \theta - \cos^2 \theta ) \cos R \phi,  \label{RHwave1}\\
v &= -  {K} R \cos^{R-1} \theta \sin \theta \sin R \phi  ,   \label{RHwave2}
\end{align}
where $ {\omega} =  {K} = 7.848 \times 10^{-6}$. For $ {h}_0 = 8 \times 10^3$, the normalized $ {\eta}$ is given by
\begin{equation}
 {\eta} = \frac{1}{ {\tilde{g}}} \left [   {A}(\theta) +  {B}(\theta) \cos R \phi +   {C} (\theta) \cos 2 R \phi \right ] ,  \label{RHwave3}
\end{equation}
where $\tilde{g}$ is the normalized gravitational constant and the coefficients $ {A}(\theta)$, $ {B}(\theta)$ and $ {C}(\theta)$ are defined as
\begin{align*}
 {A}(\theta) &= \frac{ {\omega}}{2} ( 2  {\Omega} +  {\omega} ) \cos^2 \theta  + \frac{ {K}^2}{4}  \cos^{2(R-1)} \theta [ ( R + 1) \cos^4 \theta + C_A \cos^2 \theta - 2 R^2 ], \\
 {B} (\theta) &= \frac{2 (  {\Omega} +  {\omega} )  {K} }{ (R+1) (R+2) } \cos^R \theta [ C_B - ( R + 1)^2 \cos^2 \theta ], \\
 {C} (\theta) &= \frac{{K}^2}{4}   \cos^{2R} \theta [ (R+1) \cos^2 \theta - (R+2) ],
\end{align*}
where $C_A$=$2 R^2 - R - 2 $ and $C_B $=$R^2 + 2 R + 2$. The wave number $R$ is chosen as $4$ and the Coriolis parameter $f$ is given by $ f = 2  {\Omega} \sin \theta$. For the disturbed Rossby--Hauritz wave that is proposed by Smith and Dritschel \cite{Smith}, the perturbed height $\eta^d$ in the equation reads
\begin{equation}
\eta^d = \eta \left (  1. 0 + \frac{  x x_0 + y y_0 + z z_0  }{40} \right ), 
\end{equation}
where $(x_0,~ y_0,~z_0) = (\cos \phi_0 \cos \theta_0,~ \sin \phi_0 \cos \theta_0, ~ \sin \theta_0)$ for $\phi_0 = 40 \pi / 180,~ \theta_0 = 50 \pi /  180$. Fig. \ref{SWERossbyFinal} displays the $\phi$-$\theta$ map of the Rossby-Haurwitz and disturbed Rossby--Haurwitz flow at $15$ days with the same size of time step size. Fig. \ref{SWERossbyFinal} confirms the validation of the scheme in Nekmesh in comparison with the previous results from different schemes \cite{Legat} \cite{Smith}. For the time marching, we set $\Delta t = 0.0001$.

\begin{figure}[ht]
  \centering
  \subfloat[Mass error]{\label{RossbyMass} \includegraphics[
  width=5cm]{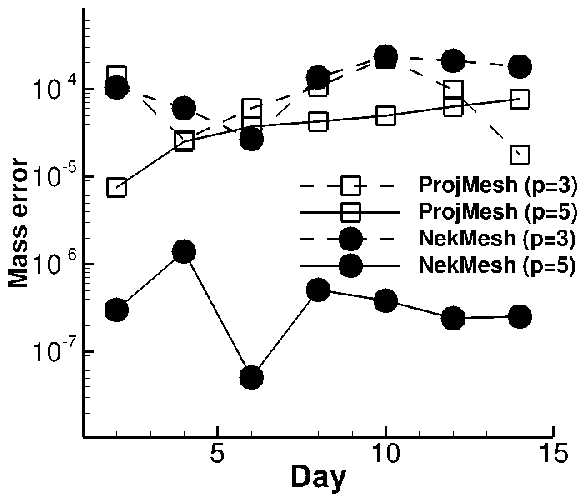}}
  \subfloat[Energy error]{\label{RossbyEnergy} \includegraphics[
  width=5cm]{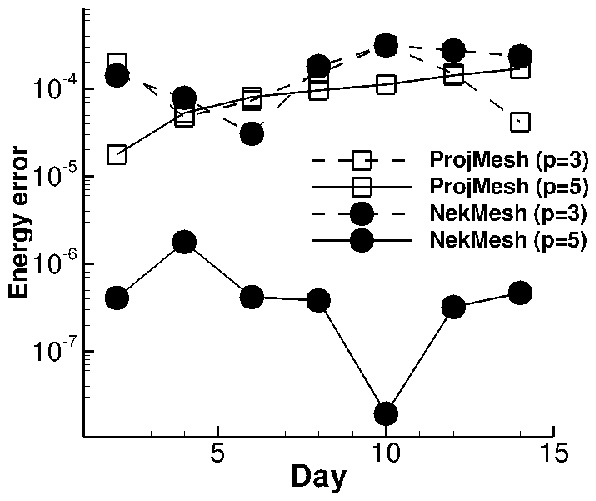}}
    \caption{Rossby--Haurwitz flow up to 14 days. $h=0.4$. }
    \label{SWERossbyTest} 
    \end{figure}

\begin{figure}[ht]
  \centering
 \includegraphics[
 width=6cm]{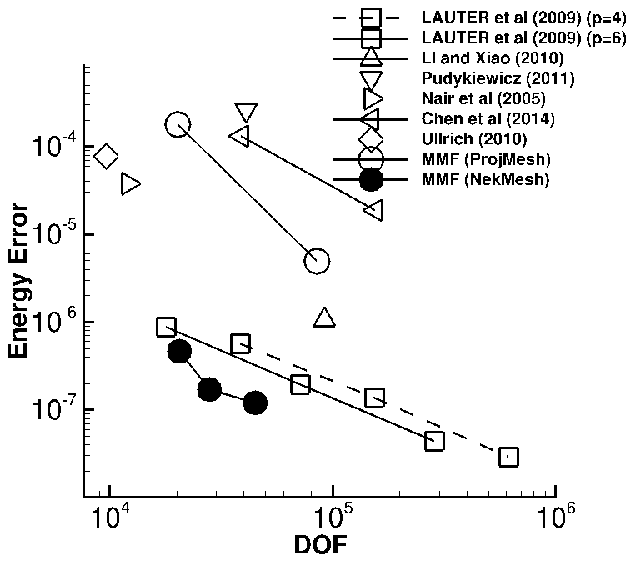}
    \caption{Comparisons of Rossby-Haurwitz flow up to 14 days. Energy error versus DOF. }
    \label{SWERossbyComparison} 
\end{figure}

In Fig. \ref{SWERossbyTest} and Table \ref{table::SWERossby}, two peculiar behaviors are observed in the convergence of NekMesh in conservational properties: First, for $p=3$, there is no significant difference in the mass and energy conservation between ProjMesh and NekMesh. We presume that the geometric approximation error is far below the discretization error, thus the discretization error dominates the mass and energy conservational error. For $p=5$, the discretization error is equivalent or smaller than the geometric approximation error, thus the geometric approximation error contributes to the mass and energy conservational error at a significant portion. At the final time $T=14.0$, the mass and energy conservation error of NekMesh is $3.3~\%$ and $2.7~\%$, respectively, compared to those of ProjMesh. Moreover, Fig. \ref{SWERossbyComparison} demonstrates that the MMF-SWE scheme on Nekmesh exhibits the most accurate energy conservation error compared with other previously proposed schemes \cite{Lauter2008,Nair2005,Li2010,Janusz2011,Chen2014,Ullrich2010}.

    \begin{table}[ht]
\resizebox{\textwidth}{!}{
\begin{tabular}{c c c c c c c c c}
 \hline\noalign{\smallskip} 
Time &   \multicolumn{4}{c}{Mass error}  &  \multicolumn{4}{c}{Energy error} \\
       &   \multicolumn{2}{c}{$p=3$}  & \multicolumn{2}{c}{$p=5$}  &    \multicolumn{2}{c}{$p=3$}  & \multicolumn{2}{c}{$p=5$}   \\
 \cmidrule(lr){2-3} \cmidrule(lr){4-5} \cmidrule(lr){6-7} \cmidrule(lr){8-9}
 & ProjMesh & NekMesh & ProjMesh & NekMesh & ProjMesh & NekMesh & ProjMesh & NekMesh \\
  \noalign{\smallskip}\hline\noalign{\smallskip}
 2 &  1.41E-04  &  1.04E-04  &  7.58E-06  &  2.98E-07  &  1.96E-04  & 1.43E-04   & 1.77E-05  &  4.10E-07   \\
 4 & 2.64E-05  & 6.01E-05   & 2.54E-05   &  1.40E-06  &  4.73E-05  &  7.84E-05  & 5.30E-05  &  1.77E-06   \\
 6 & 6.09E-05 & 2.69E-05   & 3.75E-05   & 5.04E-08   & 7.37E-05   &  3.04E-05  & 8.04E-05  &  4.20E-07   \\
 8 & 1.07E-04 & 1.36E-04   & 4.28E-05   &  5.04E-07  &  1.49E-04  &  1.80E-04  & 9.51E-05  &  3.79E-07   \\
 10 & 2.24E-04 & 2.35E-04   &  4.94E-05  &  3.82E-07  &  3.12E-04  &  3.16E-04  & 1.12E-05  & 1.93E-08   \\
  12 & 9.80E-05 & 2.90E-04   &  6.35E-05  & 2.42E-07   &  1.50E-04  &  2.77E-04  & 1.42E-04  & 3.20E-07   \\
 14 & 1.76E-05 & 1.82E-04   & 7.63E-05   &  2.55E-07  &  4.21E-05  &  2.38E-04  & 1.70E-04  &  4.65E-07  \\
\end{tabular}}
\caption{Test problem for unsteady zonal flow at 0.5 days. $h$=$0.4$ }
\label{table::SWERossby}
\end{table}

\subsection{Isolated Mountain and Unstable jet}
As the fifth Williamson's test case, which was first introduced by Takacs \cite{Takacs}, this problem solves the zonal flow over an isolated cone-shaped mountain. For the following initial velocity of zonal flow as $ \mathbf{u} = 20 / r_a \cos \theta \boldsymbol{\phi}$, the initial $ {\eta}$ is given as
\begin{equation*}
 {\eta}  =   - \frac{1}{ {\tilde{g}}} \left (  {\Omega}  {u}_0 + \frac{1}{2}  {u}_0^2 \right ) \sin^2 \theta,
\end{equation*}
and the initial $ {H}_0$ is defined as follows.
\begin{equation*}
 {H}_0 =  \frac{5960}{r_a} -   \frac{2000}{r_a}  \left ( 1 - \frac{r}{R} \right )   ,
\end{equation*}
where $R$ is chosen to be $\pi/9 $. The distance function $r$ centered at $\phi_c = 3 \pi / 2$ and $\theta_c = \pi / 6$ is obtained as follows: $r = \min \left [ R^2, (\phi - \phi_c)^2 + (\theta - \theta_c)^2 \right ]$. For the time marching, we set $\Delta t = 0.0001$.

\begin{figure}[ht]
  \centering
      \subfloat[$p=5$]{\label{ IsoMp5 } \includegraphics[
      width=5cm]{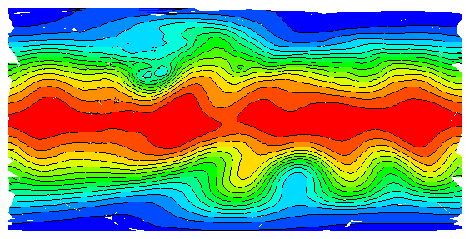}}
  \subfloat[$p=7$]{\label{IsoMp8} \includegraphics[
  width=5cm]{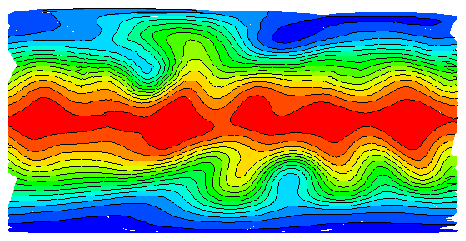}}
    \caption{Surface elevation of the isolated mountain problem at 15 days. $h$=$0.4$. Insufficient resolution at $p=5$ (left) can be solved at $p=7$ (right) without being compromised by geometric approximation error.}
    \label{SWETestISM} 
    \end{figure}

\begin{figure}[ht]
  \centering
   \includegraphics[
   width=6cm]{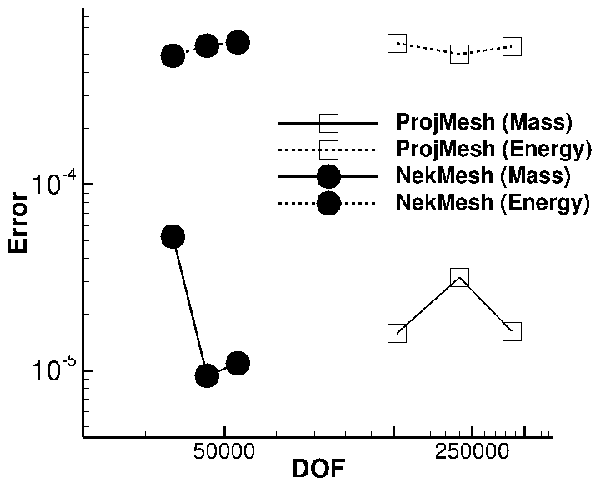}
    \caption{Mass and energy conservation error of ProjMesh and NekMesh for isolated mountain test case. }
    \label{SWEISMErr} 
    \end{figure}

A specific pattern of the surface elevation at 15 days is regenerated in the Nekmesh of $h$=$0.4$ both for $p$=$5$ and $p$=$7$ as shown in Fig \ref{SWETestISM}. For $p$=$5$, some part of the contour of surface elevation is not distinctive without shape, but for $p$=$7$ every contour is distinctive and coincident with the known distribution. The same result can be obtained by ProjMesh \cite{MMF3}. However, the total number of grid points required for NekMesh is approximately $11 \%$ of the required grid points of ProjMesh as shown in Fig. \ref{SWEISMErr}; $3.5 \times 10^4$ grid points for NekMesh and $3.2 \times 10^5$ grid points for ProjMesh. NekMesh can provide improved accuracy and less conservational loss with much less resolution even in the presence of non-differential objects such as an isolated coned mountains.

\begin{figure}[ht]
  \centering
   \includegraphics[
   width=5cm]{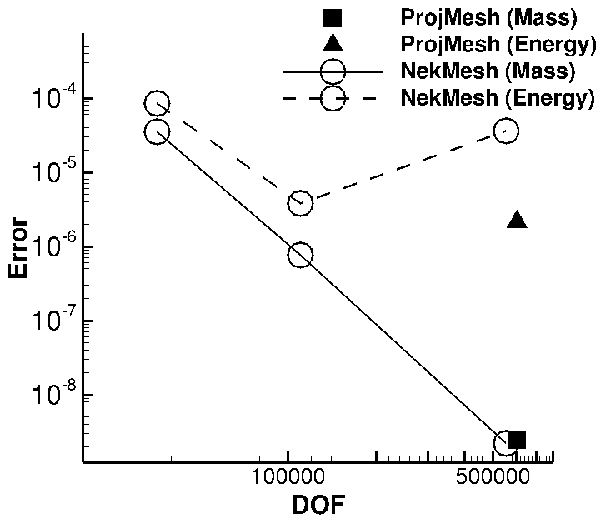}
    \caption{Mass and energy conservation error of ProjMesh and NekMesh for stable jet flow without disturbance. }
    \label{SWEJetErr} 
    \end{figure}

First proposed by Galewsky et. al. \cite{Galewsky}, the unstable jet simulation generates an unstable jet by adding a perturbation to a balanced jet. In ref. \cite{MMF3}, the unstable jet was successfully generated on the mesh of $h=0.08$ and $h=0.15$ with $p=6$ by ProjMesh. However, as observed in Fig. \ref{SWEJetErr}, the impact of geometric approximation error looks trivial if dof is sufficiently large which is closer to $5 \times 10^5$ in this test problem, or equivalently $h$ is less than $0.1$.

\begin{figure}[ht]
  \centering
      \subfloat[$h=0.08$]{\label{Jeth008 } \includegraphics[
      width=12cm]{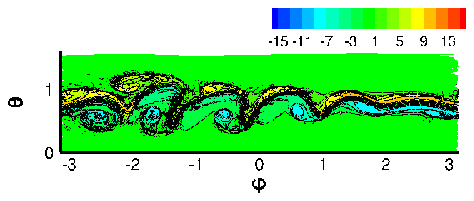}} \\
  \subfloat[$h=0.2$]{\label{Jeth02} \includegraphics[ 
  width=12cm]{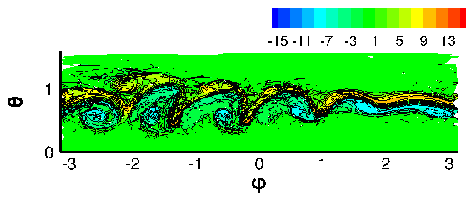}}
    \caption{Relative vorticity of the unstable jet test after 6 days. $p=6$.}
    \label{SWETestJet} 
    \end{figure}
    
\section{Discussion}
Geometric approximation error is one of the most unique phenomena in the domain with curved boundaries or on a curved surface. There have been diverse conjectures concerning the contribution of geometric approximation error to the deterioration of the accuracy in the numerical simulations, but the results of this paper seem to validate those conjectures as described as follows: 
\begin{enumerate}
\item For insufficient resolutions in which discretization errors dominate, geometric approximation error does not make any significant impact on the numerical schemes.
\item For sufficient resolutions in which discretization errors are equivalent to, or less than geometric approximation error, the presence of geometric approximation error prevents further convergence of the overall error, called the \textit{geometric approximation error saturation in the overall error}. 
\item Geometric approximation error significantly contributes to the conservation of mass and energy because the changes in the conserved quantity is approximately in the similar size as geometric approximation error. 
\item Even though the overall error is relatively larger than geometric approximation error, geometric approximation error can have a significant impact on the overall accuracy after a long time integration.
\end{enumerate}

With the typical mesh, ensuring negligible geometric approximation error involves using a very fine mesh with a very small edge length. Increasing the polynomial order does not decrease geometric approximation error. Consequently, any computational simulation, particular a long time simulation, on the domain with curved boundary or on curved domain is essentially very expensive and time-consuming. The other option as suggested in this paper is to use advanced mesh techniques such as high-order curvilinear mesh. Even with significantly small number of grid points on any curved domain, some \textit{heavy} computations can be computed with significantly less computational time and memory requirements because of the negligible effects of geometric approximation error.

\section*{Acknowledgements}
This research was supported by the Basic Science Research Program through the National Research Foundation of Korea (NRF) and funded by the Ministry of Education, Science and Technology (No. 2016R1D1A1A02937255).

\bibliographystyle{elsarticle-num}
\bibliography{MMFNekMesh_JSC}

\end{document}